\documentclass[10pt]{article}
\usepackage[cmyk]{xcolor}
\usepackage{booktabs}
\usepackage{subfigure}
\usepackage{cite}
\usepackage{graphicx}
\usepackage{multicol} %实现同一页中实现不同的分栏
\usepackage{amsmath, amsfonts, mathtools, amssymb}
\usepackage{amsthm}
\usepackage{hyperref}
\usepackage{thmtools}
\usepackage{enumerate}
\usepackage{mathrsfs}
\usepackage{cancel}
\usepackage{lmodern,txfonts}
\usepackage{type1cm}
\usepackage{bm}
\usepackage{lipsum}
\usepackage{setspace}
\usepackage{geometry}
\makeatletter

\newcommand{\Rmnum}[1]{\expandafter\@slowromancap\romannumeral #1@}
\makeatother

\newtheorem{remark}{Remark}[section]

\numberwithin{equation}{section}

\usepackage{color,bm} \usepackage[normalem]{ulem}% The line notes delete

\usepackage[T1]{fontenc}
\usepackage[utf8]{inputenc}
\usepackage{authblk}
\title{Asymptotic-preserving neural networks for multiscale Vlasov-Poisson-Fokker-Planck system in the high-field regime}
\author[1,2,3]{Shi Jin}
\author[1,2,4,5]{Zheng Ma}
\author[1,2]{Tian-ai Zhang \thanks{Corresponding author:zhangtianai@sjtu.edu.cn}}
\affil[1]{School of Mathematical Sciences, Shanghai Jiao Tong University, Shanghai, 200240, China}
\affil[2]{Institute of Natural Sciences, MOE-LSC, Shanghai Jiao Tong University, Shanghai, 200240, China}
\affil[3]{Shanghai Artificial Intelligence Laboratory, Shanghai, China}
\affil[4]{Qing Yuan Research Institute, Shanghai Jiao Tong University, Shanghai, 200240, China}
\affil[5]{CMA-Shanghai, Shanghai Jiao Tong University, Shanghai, China}

\date{}
\begin{document}
\maketitle
%%%%%%%%%%%%%%%%%%%%%%%%%%%%%%%%%%%%%%%%%%%%%%%%
%%%%%%%%%%%%%%%%%%%%%%%%%%%%%%%%%%%%%%%%%%%%%%%%
\begin{abstract}
The Vlasov-Poisson-Fokker-Planck (VPFP) system is a fundamental model in plasma physics that describes the Brownian motion of a large ensemble of particles within a surrounding bath. Under the high-field scaling, both collision and field are dominant. This paper introduces two Asymptotic-Preserving Neural Network (APNN) methods within a physics-informed neural network (PINN) framework for solving the VPFP system in the high-field regime. These methods aim to overcome the computational challenges posed by high dimensionality and multiple scales of the system. The first APNN method leverages the micro-macro decomposition model of the original VPFP system, while the second is based on the mass conservation law. Both methods ensure that the loss function of the neural networks transitions naturally from the kinetic model to the high-field limit model, thereby preserving the correct asymptotic behavior. Through extensive numerical experiments, these APNN methods demonstrate their effectiveness in solving multiscale and high dimensional uncertain problems, as well as their broader applicability for problems with long time duration and non-equilibrium initial data.
\end{abstract}

%%%%%%%%%%%%%%%%%%%%%%%%%%%%%%%%%%%%%%%%%%%%%%%%
%%%%%%%%%%%%%%%%%%%%%%%%%%%%%%%%%%%%%%%%%%%%%%%%
\section{Introduction}

The Vlasov-Poisson-Fokker-Planck (VPFP) system is a fundamental kinetic model in plasma physics, which takes into account the interactions between electrons and a surrounding bath via the Coulomb force \cite{chandrasekhar1943stochastic}. It consists of a Liouville equation with a Fokker-Planck operator in velocity space, along with a Poisson equation corresponding to the electrostatic force \cite{victory1990classical}. The high-field limit of the kinetic equations characterizes a usual scenario where the external field, such as an electric field, is strong enough to balance the collision term \cite{jin2010asymptotic}. Section 2 provides a detailed exposition of the VPFP system and its high-field limit.

The high dimensionality of the phase space in the kenitic equations ($6$ dimensions plus time for the general cases) raises a significant challenge in simulating the VPFP system. Moreover, the uncertainties of the kinetic models can introduce more dimensions\cite{jin2017uncertainty, poette2019gpc, poette2022numerical, jin2017uniform}. The uncertainties may arise from the collision kernels, scattering coefficients, initial or boundary data, source or forcing terms and so on, which are important for the application of the kinetic models to practical systems\cite{bird1994molecular, berman1986collision, koura1991variable, jin2022asymptotic}. Many efforts have been devoted to tackling the curse of dimensionality. Classical techniques, such as the Monte Carlo methods\cite{longo2000monte, havlak1996numerical}, are hindered by their low-order accuracy or slow convergence rates. Recently, machine learning, especially deep neural networks (DNNs), shows potential in resolving the Partial Differential Equations (PDEs) of high dimensionality \cite{lagaris1998artificial, long2018pde, khoo2021solving, lu2021learning, kovachki2021neural}. The key idea is to approximate the objective functions by DNNs, and optimize the network parameters by minimizing the loss function encoding the PDEs and initial-boundary conditions. Many work has been devoted to building suitable loss functions to develope DNN methods. The popular methods are the Physics-Informed Neural Networks (PINNs) \cite{raissi2019physics} and the Deep Galerkin Method (DGM) \cite{sirignano2018dgm}, which minimize the $L^2$-residual of the PDEs and the initial-boundary conditions. The Deep Ritz Method \cite{yu2018deep} was designed to solve PDEs with variational structures by leveraging its Ritz formulation. Moreover, the new area of operator learning has risen, which involves developing novel learning frameworks to approximate the solution operator of PDEs.  Such techniques include deep operator networks (DeepONets) \cite{chen1995universal, lu2021learning}, neural operators \cite{kovachki2021neural, li2020fourier} and associated variants tailored for physics problems \cite{wang2021learning, li2021physics}. For supplementary information on machine learning methods for solving PDEs, interested readers are advised to consult the comprehensive review \cite{weinan2021algorithms}. Applications of these methods to the kinetic equations are detailed in references \cite{chen2021solving, delgadillo2021multiscale, li2023learning, lou2021physics, li2021physics}.

Besides the curse of dimensionality, the multiscale nature of the VPFP system with high-field limit poses another challenge. General DNN methods, like PINNs, have difficulties in dealing with multiscale problems \cite{jin2023asymptotic}. Specifically, PINNs only capture the leading order term related to a single scale in loss functions when the scale parameter is small \cite{jin2013asymptotic}, thus fail to capture the correct asymptotic limit equation of the VPFP system in the high-field regime as shown in Section 3. Asymptotic-Preserving (AP) schemes are effective numerical schemes for multiscale problems. AP schemes preserve the continuous asymptotic limit in a numerically uniformly stable way \cite{jin2022asymptotic}. In recent years, AP schemes are widely adopted for kinetic and hyperbolic equations with multiple time and space scales \cite{jin2010asymptotic}, and also applied in the plasma \cite{Degond2013} and the high-field regime \cite{jin2011asymptotic, crouseilles2011asymptotic, zhu2017vlasov, carrillo2021variational}. However, the classical AP schemes still encounter challenges in the high dimensional cases.

To address both the curse of dimensionality and multiscale difficulty, the primary goal of this work is to develop DNN methods with AP property, which are applied to VPFP system in the high-field regime. We propose two Asymptotic-Preserving Neural Network (APNN) methods for the VPFP system. The first is based on the micro-macro model of the VPFP system \cite{crouseilles2011asymptotic}. By reformulating the loss function of the vanilla PINN method based on the micro-macro decomposition, we construct an associated loss function with AP property. This method has been successfully applied to the linear transport equations \cite{jin2013asymptotic}, the steady radiative transfer equations \cite{lu2022solving} and so on\cite{li2022model,jin2023asymptotic}. However, this method necessitates initial data close to the local Maxwellian, since its foundation is the micro-macro decomposition. Since kinetic models are more useful when the  system is far from the equilibrium inital data \cite{jin2011asymptotic}, we try to retain this advantage. Hence, we propose the second APNN method, which is free from the explicit expression of the local Maxwellian and the requirement for the equilibrium initial data. This method is based on the mass conservation law, and has the potential for  a  broader applicability for the cases with long time duration or non-equilibrium initial data.

Although our numerical experiments are conducted in one-dimensional physical and velocity spaces, we do include the examples that involve high-dimensional uncertainties to demonstrate the ability of our APNNS for high-dimensional problems.

The paper is organized as follows: Section 2 introduces the VPFP system and its high-field limit, along with the micro-macro model of the VPFP system. Section 3 presents our two proposed APNN methods and demonstrates their AP property. In Section 4, we evaluate and compare the performance of the two proposed APNN methods through a series of numerical examples. We conclude the paper in Section 5.

%%%%%%%%%%%%%%%%%%%%%%%%%%%%%%%%%%%%%%%%%%%%%%%%
%%%%%%%%%%%%%%%%%%%%%%%%%%%%%%%%%%%%%%%%%%%%%%%%
\section{The VPFP system and its high-field limit}

This section introduces the Vlasov-Poisson-Fokker-Planck (VPFP) system and its high-field limit. The micro-macro model for the VPFP system is also presented,  which provides a suitable framework for deriving the high-field limit.

%%%%%%%%%%%%%%%%%%%% VPFP %%%%%%%%%%%%%%%%%%%%%
\subsection{The VPFP system}
In kinetic equations, a strong external field, such as the electric or gravitational field, often exists and counteracts the collision term, leading to the high-field limit \cite{cercignani1997high}. In electrostatic plasmas for example, the electrons interact with the surrounding bath through the Coulomb force. Thus the time evolution of the electron distribution function $f:(t, x, v) \in \mathbb{R}_{+} \times \mathbb{R}^N \times \mathbb{R}^N \rightarrow \mathbb{R}_{+}$ is governed by the VPFP equations, with the action of a self-consistent potential $\phi(t,x)$. Specifically, the VPFP system is:
\begin{subequations}
    \begin{equation}\label{vlasov}
        \partial_t f+v\cdot\nabla_x f-\frac{1}{\varepsilon}\nabla_x \phi \cdot \nabla_v f =\frac{1}{\varepsilon} \nabla_{v} \cdot\left[v f+\nabla_{v} f\right] := \frac{1}{\varepsilon} \mathscr{Q}(f),
    \end{equation}
    \begin{equation}\label{poisson}
        -\triangle_x \phi(t, x) =\rho(t, x)-h(x),
    \end{equation}
\end{subequations}
where 
\begin{equation}\label{rho-def}
\rho(t, x)=\int_{\mathbb{R}^N} f(t, x, v) \mathrm{d} v
\end{equation}
is the density of electrons. $\mathscr{Q}$ is defined as linear Fokker-Planck operator:
\begin{equation}\label{fpop}
    \mathscr{Q}(f)(t, x, v)=\nabla_v\cdot\left[v f(t, x, v)+\nabla_v f(t, x, v)\right].
\end{equation}
The function $h(x)$ represents a given positive background charge density. Hence the assumed global neutrality relation is:
\begin{equation}
    \int_{\mathbb{R}^N} \int_{\mathbb{R}^N} f^0(x, v) \mathrm{d} x \mathrm{~d} v=\int_{\mathbb{R}^N} h(x) \mathrm{d} x .
\end{equation}
The parameter $\varepsilon$ denotes the ratio between the mean free path and the Debye length. The limiting process $\varepsilon \rightarrow 0$ indicates the high-field limit of the VPFP system. In the high-field regime, the strong forcing term balances the Fokker-Planck diffusion term \cite{arnold2001low}, which is different from the low-field limit (or named as parabolic limit) \cite{poupaud2000parabolic}.

%%%%%%%%%%%%%%%%%%%% High-field %%%%%%%%%%%%%%%%%%%%%
\subsection{The high-field limit}
In this section, we give the limit equation of the VPFP system \ref{vlasov}-\ref{poisson} at $\varepsilon \rightarrow 0$. First, integrating the Vlasov equation \ref{vlasov} over $v$ in $\mathbb{R}^N$ gives:
\begin{equation}
    \partial_t \int_{\mathbb{R}^N} f d v+\nabla_{x} \cdot \int_{\mathbb{R}^N} v f d v-\frac{1}{\varepsilon} \int_{\mathbb{R}^N} \nabla_v \cdot\left(\nabla_{x} \phi f d v\right)=\frac{1}{\varepsilon} \int_{\mathbb{R}^N} \nabla_{v} \cdot\left(v f+\nabla_{v} f\right) d v.
\end{equation}
Integrating by parts, one has:
\begin{equation}\label{rhotrans}
    \partial_t \rho+\nabla_x \cdot j=0,
\end{equation}
where the flux $j$ is defined as:
\begin{equation}
    j=\int_{\mathbb{R}^N} v f(t, x, v) \mathrm{d} v.
\end{equation}
Next, one multiplies the Vlasov equation \eqref{vlasov} by $v$ then integrates over $\mathbb{R}^N$, and takes the limit $\varepsilon \rightarrow 0$, yielding:
\begin{equation}
    0=\int_{\mathbb{R}^N} f \nabla_{x} \phi+v f+\nabla_v f d v.
\end{equation}
Deriving this equation gives:
\begin{equation}
    j=-\rho\left(\nabla_{\mathbf{x}} \phi\right).
\end{equation}
Substituting it in  equation \eqref{rhotrans} yields the high-field limit equation:
\begin{equation}\label{limit}
    \left\{\begin{array}{l}
    \partial_t \rho-\nabla_x \cdot\left(\rho \nabla_x \phi\right)=0,\\
    -\triangle_x \phi=\rho-h(x).
    \end{array}\right.
\end{equation}
The rigorous proof for the high-field limit of the VPFP system in one-dimension can be found in \cite{goudon2005multidimensional, nieto2001high}.

%%%%%%%%%%%%%%%%%%%% MM %%%%%%%%%%%%%%%%%%%%%
\subsection{The micro-macro model and the high-field limit}

This section presents the micro-macro model, from which the high-field limit can be readily derived. First, an equivalent formulation of the VPFP system \eqref{vlasov}-\eqref{poisson} is:
\begin{equation}\label{newvfp}
    \begin{aligned}
        \partial_t f+v \cdot \nabla_x f &=\frac{1}{\varepsilon} \nabla_v\cdot\left[(v+\nabla_x\phi) f+\nabla_v f\right] := \frac{1}{\varepsilon} \mathscr{L} f,\\
        -\Delta_x \phi &=\rho-h(x).
    \end{aligned}
\end{equation} 
The linear operator $\mathscr{L}$ retains the formulation of a Fokker-Planck operator dependent on $\nabla_x \phi$.  The null space of $\mathscr{L}$ is:
\begin{equation}
    \mathscr{N}(\mathscr{L})=\operatorname{Span}\{\mathscr{M}\}=\{f=\rho \mathscr{M}, \text{where} \, \rho:=\langle f\rangle\},
\end{equation}
where $\mathscr{M}$ is the shifted Maxwellian depending on $(t,x)$ through the potential $\phi$ in the form:
\begin{equation}
    \mathscr{M}(t, x, v)=\frac{1}{\sqrt{2 \pi}} \exp\left(-\frac{|v+\nabla_x\phi(t, x)|^2}{2}\right).
\end{equation}
Notation $\langle f\rangle$ is defined as:
\begin{equation}
    \langle f\rangle=\int_{\mathbb{R}^N} f(v) d v.
\end{equation}
The rank of $\mathscr{L}$ is:
\begin{equation}
    \mathscr{R}(\mathscr{L})=(\mathscr{N})^{\perp}(\mathscr{L})=\{f\, \text{such that} \, \langle f\rangle=\int_{\mathbb{R}^N} f d v=0\}.
\end{equation}
The idea for the micro-macro model begins with the decomposition of $f$ into the equilibrium part and the non-equilibrium part \cite{crouseilles2011asymptotic, lemou2008new}:
\begin{equation}\label{mm}
    f=\rho \mathscr{M}+\varepsilon g,
\end{equation}
where $\rho \mathscr{M}$ is the equilibrium part of $f$. The non-equilibrium part $g$ satisfies $\langle g \rangle = 0$, which implies $g \in \mathscr{R}(\mathscr{L})$. The orthogonal projector $\Pi$ in $L^2\left(\mathscr{M}^{-1} d v\right)$ onto $\mathscr{N}(\mathscr{L})$ is defined by:
\begin{equation}\label{projection}
    \Pi \varphi=\langle\varphi\rangle \mathscr{M}.
\end{equation}

In deriving the micro-macro model, we first consider  the one-dimensional VPFP system formulated as:
\begin{equation}\label{1d}
    \partial_t f+v \partial_x f=\frac{1}{\varepsilon} \partial_v\left[(v+\partial_x\phi) f+\partial_v f\right]=: \frac{1}{\varepsilon} \mathscr{L} f.
\end{equation}
Applying the micro-macro decomposition \eqref{mm} to the equation \eqref{1d} yields:
\begin{equation}\label{mm1}
    \partial_t(\rho \mathscr{M})+\varepsilon\partial_t g+ v \partial_x(\rho \mathscr{M})+\varepsilon v \partial_x g=\mathscr{L} g.
\end{equation}
Operating $(I-\Pi)$ on  equation \eqref{mm1} gives:
\begin{equation*}
    (I-\Pi)\left[\partial_t(\rho \mathscr{M})+v \partial_x(\rho \mathscr{M})\right]+\varepsilon(I-\Pi)\left[\partial_t g+v \partial_x g\right]=(I-\Pi) \mathscr{L} g.
\end{equation*}
Since $\Pi(g)=\Pi\left(\partial_t g\right)=\Pi(\mathscr{L} g)=0$ \cite{crouseilles2011asymptotic}, and
\begin{equation*}
    \Pi\left[\partial_t(\rho \mathscr{M})\right]=\Pi\left[\mathscr{M} \partial_t \rho-\rho\left(\partial_{xt} \phi\right)(v+\partial_x\phi) \mathscr{M}\right]=\mathscr{M} \partial_t \rho,
\end{equation*}
one gets
\begin{equation}\label{mm2}
    \varepsilon\partial_t g+\varepsilon(I-\Pi)\left(v \partial_x g\right)=\mathscr{L} g-(I-\Pi)\left(v \partial_x(\rho \mathscr{M})\right)-\rho(\partial_t\mathscr{M}).
\end{equation}
Then the $\Pi$ projection \eqref{projection} is applied to the equation \eqref{mm1}
\begin{equation}
    \partial_t \rho+\partial_x\langle v (\rho \mathscr{M})\rangle+\varepsilon\partial_x\langle v g\rangle =0.
\end{equation}
Therefore, one reformulates the original VPFP system \eqref{vlasov}-\eqref{poisson} as:
\begin{equation}\label{mm-vpfp}
    \left\{\begin{aligned}
    \partial_t g+\left(v \partial_x g-\partial_x\langle v g\rangle \mathscr{M}\right) &=\frac{1}{\varepsilon}\{\mathscr{L} g- [v \partial_x(\rho \mathscr{M})\\
    &- \partial_x\langle v(\rho \mathscr{M})\rangle \mathscr{M}]-\rho(\partial_t\mathscr{M})\}, \\
    \partial_t \rho+\partial_x\langle v (\rho \mathscr{M})\rangle+\varepsilon\partial_x\langle v g\rangle &=0, \\
    -\partial_x^2 \phi &=\rho-h, \\
    \langle g\rangle &= 0.
    \end{aligned}\right.
\end{equation}
The equivalence between this model and the original VPFP system \eqref{vlasov}-\eqref{poisson} has been demonstrated in \cite{crouseilles2011asymptotic}. Note that the conservation condition $\langle g\rangle=0$ is usually satisfied by classical mesh-based AP schemes,  but it fails to hold  by DNN-predicated methods, unless incorporated into the loss function, as done in \cite{jin2023asymptotic}.

Now, the formal derivation of the high-field limit model can also be done easily from the micro-macro model \eqref{mm-vpfp}. Taking the limit as $\varepsilon \rightarrow 0$ in the second equation of the micro-macro model \eqref{mm-vpfp} yields:
\begin{equation}\label{middle}
    \partial_t \rho+\partial_x\langle v(\rho \mathscr{M})\rangle=0.
\end{equation}
Substituting the explicit expression for $\mathscr{M}$ in the equation \eqref{middle} reformulates it as:
\begin{equation}
    \partial_t \rho-\partial_x(\rho\partial_x\phi)=0.
\end{equation}
This results in the limit model:
\begin{equation}\label{limit1}
    \left\{\begin{aligned}
    \partial_t \rho-\partial_x(\rho\partial_x\phi) &=0, \\
    -\partial_x^2\phi &=\rho-h.
    \end{aligned}\right.
\end{equation}
It is exactly the high-field limit model \eqref{limit} in one-dimension, which justifies the AP property of the micro-macro model.

Finally,  as shown in \cite{crouseilles2011asymptotic, lemou2008new}, the Initial and periodic boundary conditions imposed on $f$ can be directly translated into boundary conditions for $\rho$ and $g$.

%%%%%%%%%%%%%%%%%%%%%%%%%%%%%%%%%%%%%%%%%%%%%%%%
%%%%%%%%%%%%%%%%%%%%%%%%%%%%%%%%%%%%%%%%%%%%%%%%
\section{Two formulations of asymptotic-preserving neural networks (APNNs)}

We seek to approximate the solution to the VPFP system by neural networks. Under the framework of PINNs and other neural network-based methods, loss functions encoding the PDE constraints are minimized to converge on solutions. PINNs is our starting point but vanilla DNN approximations face limitations in tackling multiscale problems where both microscopic and macroscopic dynamics are involved. To deal with this difficulty, asymptotic-preserving neural networks (APNNs) have been proposed \cite{jin2023asymptotic}. The key idea is to design a loss with AP property. Specifically, the loss function automatically transits from kinetic system to it asymptotic limit system when scale parameter $\varepsilon\rightarrow0$. In practice, the original multiscale system is first reformulated into an equivalent model with the AP property. This modified form is then embedded into the vanilla PINN loss function, with initial and boundary conditions treated as regularization terms. 

To improve multiscale VPFP approximations of vanilla PINNs, we propose two types of loss functions with AP property. The constructions of the novel AP loss functions are respectively inspired by the derivations of the high-field limit models in Sections 1.2 and 1.3. 

Consider the VPFP system with initial and boundary conditions over a bounded domain as an approximation of the original system \eqref{vlasov}-\eqref{poisson}, in which $(t, x, v) \in \mathcal{T} \times \mathcal{D} \times \Omega$ and $\Omega$ is symmetric in $v$:
\begin{equation}\label{ibcVPFP}
    \begin{cases}
        \partial_t f+v \cdot \nabla_x f=\frac{1}{\varepsilon} \mathscr{L} f, & (t, x, v) \in \mathcal{T} \times \mathcal{D} \times \Omega, \\ 
        -\triangle_x \phi(t, x) =\rho(t, x)-h(x), & (t, x) \in \mathcal{T} \times \mathcal{D},\\
        \mathcal{B} f=F_{\mathrm{B}}, & (t, x, v) \in \mathcal{T} \times \partial \mathcal{D} \times \Omega, \\ 
        \mathcal{B} \phi= \Phi_{\mathrm{B}}, & (t, x) \in \mathcal{T} \times \partial \mathcal{D}\\
        \mathcal{I} f=f_0, & (t, x, v) \in\{t=0\} \times \mathcal{D} \times \Omega,\\
        \mathcal{I} \phi=  \phi_0, & (t, x) \in\{t=0\} \times \mathcal{D}
    \end{cases}
\end{equation}
where $F_{\mathrm{B}}, f_0, \Phi_{\mathrm{B}}, \phi_0$ are given functions; $\partial \mathcal{D}$ is the boundary of $\mathcal{D}$; and $\mathcal{B}, \mathcal{I}$ are boundary and initial operators, respectively.  In the following sections, the limitations of vanilla PINNs for the multiscale VPFP system will be explained. We will also provide the detailed constructions of the two proposed AP loss functions and formally prove their AP properties.

%%%%%%%%%%%%%%%%%%%% PINN %%%%%%%%%%%%%%%%%%%%%
\subsection{The failure of PINNs to resolve small scales}

In this section, the deficiency of vanilla PINN are demonstrated in the case with small scale parameter $\varepsilon$. First we use DNNs to  parameterize two functions $f$, $\phi$ and denote the set of neural network parameters by $\eta$ for brevity. Then two neural networks are implemented:
\begin{equation}
    f_\gamma^{\mathrm{NN}}(t, x):= \log\left(1 + \exp \left(\tilde{f}_\gamma^{\mathrm{NN}}(t, x)\right)\right) \approx f(t, x), \,\,\, \text{and}\,\,\,\, \phi_\gamma^{\mathrm{NN}}(t, x)\approx \phi(t, x).
\end{equation}
$\tilde{f}_\gamma^{\mathrm{NN}}$ and $\phi_\gamma^{\mathrm{NN}}$ are both fully-connected neural networks. Note that the setting for $f_\gamma^{\mathrm{NN}}$ guarantees the non-negativity property of distribution function $f$. 

The vanilla PINN loss is constructed by the least square of the residuals of the original VPFP system, combined with penalty terms encoding the initial and boundary conditions:
\begin{equation}\label{loss PINN}
    \begin{aligned}
        \mathcal{R}_{\mathrm{PINN}}^{\varepsilon} &= \frac{\mu_1^\text{Vlasov}}{|\mathcal{T} \times \mathcal{D} \times \Omega|} \int_{\mathcal{T}} \int_{\mathcal{D}} \int_{\Omega}\left|\varepsilon\partial_tf_\gamma^{\mathrm{NN}}+\varepsilon v\partial_xf_\gamma^{\mathrm{NN}}-\mathscr Lf_\gamma^{\mathrm{NN}} \right|^2\mathrm{~d}v\mathrm{~d} x \mathrm{d} t\\
        &+ \frac{\mu_1^\text{Poisson}}{|\mathcal{T} \times \mathcal{D}|} \int_{\mathcal{T}}\int_{\mathcal{D}}\left|-\partial_x^2\phi_\gamma^{\mathrm{NN}}-(\rho_\gamma^{\mathrm{NN}}-h)\right|^2\mathrm{~d} x \mathrm{d} t \\
        &+\frac{\mu_2^f}{|\mathcal{T} \times \partial \mathcal{D} \times \Omega\mid} \int_{\mathcal{T}} \int_{\partial \mathcal{D}} \int_{\Omega}\left|\mathcal{B}f_\gamma^{\mathrm{NN}}-F_{\mathrm{B}}\right|^2 \mathrm{~d} v \mathrm{d} x \mathrm{d} t \\
        &+\frac{\mu_2^\phi}{|\mathcal{T} \times \partial \mathcal{D}|} \int_{\mathcal{T}} \int_{\partial \mathcal{D}}\left|\mathcal{B}\phi_\gamma^{\mathrm{NN}}-\Phi_{\mathrm{B}}\right|^2 \mathrm{d} x \mathrm{d} t \\
        & +\frac{\mu_3^f}{|\mathcal{D} \times \Omega|} \int_{\mathcal{D}} \int_{\Omega}\left|\mathcal{I}f_\gamma^{\mathrm{NN}}-f_0\right|^2 \mathrm{~d} v \mathrm{d} x +\frac{\mu_3^\phi}{|\mathcal{D}|} \int_{\mathcal{D}} \left|\mathcal{I}\phi_\gamma^{\mathrm{NN}}-\phi_0\right|^2\mathrm{d} x,
    \end{aligned}
\end{equation}
where $\mu_1$, $\mu_2$ and $\mu_3$ are the penalty weights to be tuned. For brevity, a formal notation is introduced:
\begin{equation}
    \mathcal{R}_{\mathrm{PINN}}^{\varepsilon} := \mu_1\mathcal{R}_{\text{residual}}^{\varepsilon} + \mu_2\mathcal{R}_{\text{bc}}^{\varepsilon} + \mu_3\mathcal{R}_{\text{ic}}^{\varepsilon}.
\end{equation}

We now examine whether this loss of PINN method possesses the AP property. Wee only need to focus on the first term of \eqref{loss PINN}:
\begin{equation}
    \begin{aligned}
        \mathcal{R}_{\text{residual}}^{\varepsilon}=
        & \frac{1}{|\mathcal{T} \times \mathcal{D} \times \Omega|} \int_{\mathcal{T}} \int_{\mathcal{D}} \int_{\Omega}\left|\varepsilon\partial_tf_\gamma^{\mathrm{NN}}+\varepsilon v\partial_xf_\gamma^{\mathrm{NN}}-\mathscr Lf_\gamma^{\mathrm{NN}} \right|^2\mathrm{~d}v\mathrm{~d} x \mathrm{d} t\\
        &+ \frac{1}{|\mathcal{T} \times \mathcal{D}|} \int_{\mathcal{T}}\int_{\mathcal{D}}\left|-\partial_x^2\phi_\gamma^{\mathrm{NN}}-(\rho_\gamma^{\mathrm{NN}}-h)\right|^2\mathrm{~d} x \mathrm{d} t .
    \end{aligned}
\end{equation}
Taking $\varepsilon \rightarrow 0$, formally this leads to:
\begin{equation}
    \begin{aligned}
        \mathcal{R}_{\text{residual}}^{\varepsilon}=
        &\frac{1}{|\mathcal{T} \times \mathcal{D} \times \Omega|} \int_{\mathcal{T}} \int_{\mathcal{D}} \int_{\Omega}\left|\mathscr Lf_\gamma^{\mathrm{NN}} \right|^2\mathrm{~d}v\mathrm{~d} x \mathrm{d} t\\
        &+ \frac{1}{|\mathcal{T} \times \mathcal{D}|} \int_{\mathcal{T}}\int_{\mathcal{D}}\left|-\partial_x^2\phi_\gamma^{\mathrm{NN}}-(\rho_\gamma^{\mathrm{NN}}-h)\right|^2\mathrm{~d} x \mathrm{d} t.
    \end{aligned}
\end{equation}
It is the least square loss of the following system: 
\begin{equation}
    \left\{\begin{aligned}
    &\mathscr{L} f = 0, \\
    &-\partial_x^2 \phi =\rho-h.
    \end{aligned}\right.
\end{equation}
At the small scaling of $\varepsilon$, our deviation shows that the PINN loss is dominated by its leading order terms. We are actually solving the equation $\mathscr{L} f=0$, which gives $f=\rho\mathscr{M}$. Clearly, this equation fails to produce the desired high-field limit equation \eqref{limit1}. Therefore, training neural networks with vanilla PINN loss may lead to an inaccurate solution when $\varepsilon$ is small. This will be furthermore, demonstrated by numerical experiments in Section 3.

%%%%%%%%%%%%%%%%%%%% MM-APNN %%%%%%%%%%%%%%%%%%%%%
\subsection{The micro-macro decomposition based APNN method}

We now present an APNN method based on the micro-macro decomposition for the VPFP system. The main idea is to utilize PINN to solve the micro-macro system \eqref{mm-vpfp} rather than the original system \eqref{vlasov}-\eqref{poisson}. DNNs are employed to parameterize functions $\rho$, $g$, $\phi$. Three independent neural networks are adopted for three functions. The set of parameters and the parameterized functions are separately denoted by $\theta$ and $\rho_\theta^\text{NN}, g_\theta^\text{NN},$ $\phi_\theta^\text{NN}$. In detail, $\rho(t,x)$ is parameterized as:
\begin{equation}
    \rho_\theta^{\mathrm{NN}}(t, x):= \log\left(1 + \exp \left(\tilde{\rho}_\theta^{\mathrm{NN}}(t, x)\right)\right) \approx \rho(t, x),
\end{equation}
which conserves the non-negative property of density $\rho$. And $\phi(t, x)$ is parameterized as:
\begin{equation}
    \phi_\theta^{\mathrm{NN}}(t, x)\approx \phi(t, x).
\end{equation}
Moreover, $g(t, x, v)$ is parameterized as: 
\begin{equation}\label{g}
    g_\theta^{\mathrm{NN}}(t, x, v) = \tilde{g}_\theta^{\mathrm{NN}}(t, x, v) - \frac{1}{|\Omega|}\int_{\Omega}\tilde{g}_\theta^{\mathrm{NN}}(t, x, v)dv\approx g(t, x, v).
\end{equation}
Here $\tilde{\rho}_\theta^{\mathrm{NN}}$ and $\tilde{g}_\theta^{\mathrm{NN}}$ and $\phi_\theta^{\mathrm{NN}}$ are all fully-connected neural networks. 

The design of $g_\theta^{\mathrm{NN}}$ in \eqref{g} automatically guarantees the conservation condition in the micro-macro model \eqref{mm-vpfp}:
\begin{equation}
    \left\langle g_\theta^{\mathrm{NN}}\right\rangle=0, \quad \forall t, x.
\end{equation}
As demonstrated in \cite{jin2023asymptotic}, this representation of $g_\theta^{\mathrm{NN}}$ improves the accuracy of predicted results, compared to the method of incorporating this conservation mechanism as a soft constraint in loss. To verify the conservation property of $g_\theta^{\mathrm{NN}}$, first, $\langle \tilde{g}_\theta^{\mathrm{NN}}\rangle :=\int_{\mathbb{R}}\tilde{g}_\theta^{\mathrm{NN}}dv$ is approximated by $\int_{\Omega}\tilde{g}_\theta^{\mathrm{NN}}dv$ in bounded domain. Then one has:
\begin{equation}
    \begin{aligned}
        g_\theta^{\mathrm{NN}}(t, x, v) &= \tilde{g}_\theta^{\mathrm{NN}}(t, x, v) - \frac{1}{|\Omega|}\int_{\Omega}\tilde{g}_\theta^{\mathrm{NN}}(t, x, v)dv\\
        &\approx \tilde{g}_\theta^{\mathrm{NN}}(t, x, v)-\frac{1}{|\Omega|}\langle \tilde{g}_\theta^{\mathrm{NN}}(t, x, v)\rangle.
    \end{aligned}
\end{equation}
Finally, the following can be obtained:
\begin{equation}
    \langle g_\theta^{\mathrm{NN}}\rangle \approx\int_{\Omega}g_\theta^{\mathrm{NN}}dv =\int_{\Omega}\tilde{g}_\theta^{\mathrm{NN}}dv - \frac{1}{|\Omega|}(\int_{\Omega}\tilde{g}_\theta^{\mathrm{NN}}dv)(\int_\Omega dv)=0
\end{equation}
which confirms the conservation property.

To construct the loss function for the APNN method, we take the least squares of the residuals from the micro-macro system \eqref{mm-vpfp}:
\begin{equation}
    \begin{aligned}\label{MMloss}
        \mathcal{R}_{\mathrm{APNN}}^{\varepsilon}=& \frac{\lambda_1^\text{macro}}{|\mathcal{T} \times \mathcal{D}|} \int_{\mathcal{T}}\int_{\mathcal{D}}\left|\partial_t \rho_\theta^{\mathrm{NN}} + \partial_x\langle v(\rho_\theta^{\mathrm{NN}}\mathscr M_\theta^{\mathrm{NN}})\rangle +\varepsilon\partial_x\langle vg_\theta^{\mathrm{NN}}\rangle \right|^2\mathrm{~d} x \mathrm{d} t \\
        &+\frac{\lambda_1^\text{micro}}{|\mathcal{T} \times \mathcal{D} \times \Omega|} \int_{\mathcal{T}} \int_{\mathcal{D}} \int_{\Omega} \left| \varepsilon\partial_t g_\theta^{\mathrm{NN}} + \varepsilon(v\partial_x g_\theta^{\mathrm{NN}}-\partial_x\langle vg_\theta^{\mathrm{NN}}\rangle\mathscr M_\theta^{\mathrm{NN}})\right. \\
        &\left. - \{\mathscr{L} g_\theta^{\mathrm{NN}}-v \partial_x(\rho_\theta^{\mathrm{NN}} \mathscr{M}_\theta^{\mathrm{NN}})+\partial_x\langle v(\rho_\theta^{\mathrm{NN}} \mathscr{M}_\theta^{\mathrm{NN}})\rangle \mathscr{M}_\theta^{\mathrm{NN}}-\rho_\theta^{\mathrm{NN}}(\partial_t\mathscr{M}_\theta^{\mathrm{NN}})\}\right|^2\mathrm{~d} v \mathrm{d} x \mathrm{d} t \\
        &+ \frac{\lambda_1^\text{Poisson}}{|\mathcal{T} \times \mathcal{D}|} \int_{\mathcal{T}}\int_{\mathcal{D}}\left|-\partial_x^2\phi_\theta^{\mathrm{NN}}-(\rho_\theta^{\mathrm{NN}}-h)\right|^2\mathrm{~d} x \mathrm{d} t \\
        &+\frac{\lambda_2^f}{|\mathcal{T} \times \partial \mathcal{D} \times \Omega \mid} \int_{\mathcal{T}} \int_{\partial \mathcal{D}} \int_{\Omega}\left|\mathcal{B}( \rho_\theta^{\mathrm{NN}}\mathscr{M}_\theta^{\mathrm{NN}} + \varepsilon g_\theta^{\mathrm{NN}})-F_{\mathrm{B}}\right|^2 \mathrm{~d} v \mathrm{d} x \mathrm{d} t \\
        &+\frac{\lambda_2^\phi}{|\mathcal{T} \times \partial \mathcal{D}|} \int_{\mathcal{T}} \int_{\partial \mathcal{D}}\left|\mathcal{B}\phi_\theta^{\mathrm{NN}}-\Phi_{\mathrm{B}}\right|^2 \mathrm{d} x \mathrm{d} t \\
        & +\frac{\lambda_3^f}{|\mathcal{D} \times \Omega|} \int_{\mathcal{D}} \int_{\Omega}\left|\mathcal{I}(\rho_\theta^{\mathrm{NN}}\mathscr{M}_\theta^{\mathrm{NN}} + \varepsilon g_\theta^{\mathrm{NN}})-f_0\right|^2 \mathrm{~d} v \mathrm{d}x\\
        &+\frac{\lambda_3^\phi}{|\mathcal{D}|} \int_{\mathcal{D}} \int_{\Omega}\left|\mathcal{I}\phi_\theta^{\mathrm{NN}}-\phi_0\right|^2\mathrm{d}x,
    \end{aligned}
\end{equation}
where $\mathscr M_\theta^{\mathrm{NN}}$ can be calculated from $\rho_\theta^{\mathrm{NN}}$, $g_\theta^{\mathrm{NN}}$ and $\phi_\theta^{\mathrm{NN}}$. The loss is formally denoted as: 
\begin{equation}
    \mathcal{R}_{\mathrm{APNN}}^{\varepsilon} := \lambda_1\mathcal{R}_{\text{residual}}^{\varepsilon} + \lambda_2\mathcal{R}_{\text{bc}}^{\varepsilon} + \lambda_3\mathcal{R}_{\text{ic}}^{\varepsilon}.
\end{equation}    

Here, we consider the periodic BCs. As precised at the end of Section 1.3, ICs and the periodic BCs of $\rho$ and $g$ can be directly deduced from the IC and periodic BC imposed on $f$. Hence the initial and boundary conditions for $f$ in the bounded domain VPFP system \eqref{ibcVPFP} are reformulated as:
\begin{equation}
    \begin{cases} 
        \mathcal{B} \rho= P_{\mathrm{B}}, & (t, x) \in \mathcal{T} \times \partial \mathcal{D}\\ 
        \mathcal{B} g=G_{\mathrm{B}}, & (t, x, v) \in \mathcal{T} \times \partial \mathcal{D} \times \Omega \\
        \mathcal{B} \phi= \Phi_{\mathrm{B}}, & (t, x) \in \mathcal{T} \times \partial \mathcal{D}\\ 
        \mathcal{I} \rho=\rho_0, & (t, x) \in\{t=0\} \times \mathcal{D}\\
        \mathcal{I} g=g_0, & (t, x, v) \in\{t=0\} \times \mathcal{D} \times \Omega\\
        \mathcal{I} \phi=  \phi_0, & (t, x) \in\{t=0\} \times \mathcal{D}
    \end{cases}
\end{equation}
where $P_{\mathrm{B}}, G_{\mathrm{B}}, \Phi_{\mathrm{B}},\rho_0, g_0, \phi_0$ can be obtained from the corresponding settings of $f$. In this case, the loss \eqref{MMloss} is rewritten as:
\begin{equation}
    \begin{aligned}
        \mathcal{R}_{\mathrm{APNN}}^{\varepsilon} &=
        \lambda_1\mathcal{R}_{\text{residual}}^{\varepsilon} +\frac{\lambda_2^\rho}{|\mathcal{T} \times \partial \mathcal{D}|} \int_{\mathcal{T}} \int_{\partial \mathcal{D}}\left|\mathcal{B}\rho_\theta^{\mathrm{NN}}-P_{\mathrm{B}}\right|^2 \mathrm{d} x \mathrm{d} t\\
        &+\frac{\lambda_2^g}{|\mathcal{T} \times \partial \mathcal{D} \times \Omega \mid} \int_{\mathcal{T}} \int_{\partial \mathcal{D}} \int_{\Omega}\left|\mathcal{B}g_\theta^{\mathrm{NN}}-G_{\mathrm{B}}\right|^2 \mathrm{~d} v \mathrm{d} x \mathrm{d} t +\frac{\lambda_2^\phi}{|\mathcal{T} \times \partial \mathcal{D}|} \int_{\mathcal{T}} \int_{\partial \mathcal{D}}\left|\mathcal{B}\phi_\theta^{\mathrm{NN}}-\Phi_{\mathrm{B}}\right|^2 \mathrm{d} x \mathrm{d} t \\
        &+\frac{\lambda_3^\rho}{|\mathcal{D}|} \int_{\mathcal{D}} \int_{\Omega}\left|\mathcal{I}\rho_\theta^{\mathrm{NN}}-\rho_0\right|^2\mathrm{d}x 
        +\frac{\lambda_3^g}{|\mathcal{D} \times \Omega|} \int_{\mathcal{D}} \int_{\Omega}\left|\mathcal{I}g_\theta^{\mathrm{NN}}-g_0\right|^2 \mathrm{~d} v \mathrm{d}x \\
        &+\frac{\lambda_3^\phi}{|\mathcal{D}|} \int_{\mathcal{D}} \int_{\Omega}\left|\mathcal{I}\phi_\theta^{\mathrm{NN}}-\phi_0\right|^2\mathrm{d}x .
    \end{aligned}
\end{equation}

Wee now formally show the AP property of the loss function \eqref{MMloss}. Consider small $\varepsilon$, we only need to focus on the residual term:
\begin{equation}
    \begin{aligned}
        \mathcal{R}_{\text{residual}}^{\varepsilon} &= \frac{1}{|\mathcal{T} \times \mathcal{D}|} \int_{\mathcal{T}}\int_{\mathcal{D}}\left|\partial_t \rho_\theta^{\mathrm{NN}} + \partial_x\langle v(\rho_\theta^{\mathrm{NN}}\mathscr M_\theta^{\mathrm{NN}})\rangle +\varepsilon\partial_x\langle vg_\theta^{\mathrm{NN}}\rangle \right|^2\mathrm{~d} x \mathrm{d} t \\
        &+\frac{1}{|\mathcal{T} \times \mathcal{D} \times \Omega|} \int_{\mathcal{T}} \int_{\mathcal{D}} \int_{\Omega} \mid \varepsilon\partial_t g_\theta^{\mathrm{NN}} + \varepsilon(v\partial_x g_\theta^{\mathrm{NN}}-\partial_x\langle vg_\theta^{\mathrm{NN}}\rangle\mathscr M_\theta^{\mathrm{NN}}) \\
        &- \{\mathscr{L} g_\theta^{\mathrm{NN}}-v \partial_x(\rho_\theta^{\mathrm{NN}} \mathscr{M}_\theta^{\mathrm{NN}})+\partial_x\langle v(\rho_\theta^{\mathrm{NN}} \mathscr{M}_\theta^{\mathrm{NN}})\rangle \mathscr{M}_\theta^{\mathrm{NN}}-\rho_\theta^{\mathrm{NN}}(\partial_t\mathscr{M}_\theta^{\mathrm{NN}})\}|^2\mathrm{~d} v \mathrm{d} x \mathrm{d} t\\
        &+ \frac{1}{|\mathcal{T} \times \mathcal{D}|} \int_{\mathcal{T}}\int_{\mathcal{D}}\left|-\partial_x^2\phi_\theta^{\mathrm{NN}}-(\rho_\theta^{\mathrm{NN}}-h)\right|^2\mathrm{~d} x \mathrm{d} t.
    \end{aligned}
\end{equation}
Taking $\varepsilon \rightarrow 0$, it formally becomes:
\begin{equation}
    \begin{aligned}
        \mathcal{R}_{\text{residual}}^{\varepsilon} &= \frac{1}{|\mathcal{T} \times \mathcal{D}|} \int_{\mathcal{T}}\int_{\mathcal{D}}\left|\partial_t \rho_\theta^{\mathrm{NN}} + \partial_x\langle v(\rho_\theta^{\mathrm{NN}}\mathscr M_\theta^{\mathrm{NN}})\rangle\right|^2\mathrm{~d} x \mathrm{d} t \\
        &+\frac{1}{|\mathcal{T} \times \mathcal{D} \times \Omega|} \int_{\mathcal{T}} \int_{\mathcal{D}} \int_{\Omega} \mid - \{\mathscr{L} g_\theta^{\mathrm{NN}}-v \partial_x(\rho_\theta^{\mathrm{NN}} \mathscr{M}_\theta^{\mathrm{NN}})+\partial_x\langle v(\rho_\theta^{\mathrm{NN}} \mathscr{M}_\theta^{\mathrm{NN}})\rangle \mathscr{M}_\theta^{\mathrm{NN}}-\rho_\theta^{\mathrm{NN}}(\partial_t\mathscr{M}_\theta^{\mathrm{NN}})\}|^2\mathrm{~d} v \mathrm{d} x \mathrm{d} t\\
        &+ \frac{1}{|\mathcal{T} \times \mathcal{D}|} \int_{\mathcal{T}}\int_{\mathcal{D}}\left|-\partial_x^2\phi_\theta^{\mathrm{NN}}-(\rho_\theta^{\mathrm{NN}}-h)\right|^2\mathrm{~d} x \mathrm{d} t.
    \end{aligned}
\end{equation}
This equation coincides with the least square loss of the system:
\begin{equation}\label{mid high}
    \left\{\begin{aligned}
    &\mathscr{L} g = \left[v \partial_x(\rho \mathscr{M})- \partial_x\langle v(\rho \mathscr{M})\rangle \mathscr{M}\right]-\rho(\partial_t\mathscr{M}), \\
    &\partial_t \rho+\partial_x\langle v (\rho \mathscr{M})\rangle =0, \\
    &-\partial_x^2 \phi =\rho-h.
    \end{aligned}\right.
\end{equation}
Similar to Section 1.3, the second equation in \eqref{mid high} produces $\partial_t \rho-\partial_x\left(\rho \partial_x \phi\right)=0$. Combining it with $-\partial_x^2 \phi = \rho - h$, the high-field limit system \eqref{limit1} is derived. The AP property of this proposed DNN method based on the micro-macro decomposition is demonstrated. Therefore, this approach enables accurate predictions for both kinetic and high-field regimes.

%%%%%%%%%%%%%%%%%%%% MC-APNN %%%%%%%%%%%%%%%%%%%%%
\subsection{The mass conservation based APNN method}

In this section, we propose another APNN method. This method is built upon an equivalent model of VPFP system, which incorporates the inherent mass conversion law to the origin VPFP system. Our new APNN method does not depend upon the explicit expression of the local Maxwellian, and nor it requires  the initial data to be a local Maxwellian, thereby significantly broadening its potential applications.

\subsubsection{Enforcing the mass conservation law}

%In this part, the new model based on the mass conservation law and the straightforward derivation of its high-field limit are introduced. 
A distinguished feature of this model is to enforce the local mass conversion law in the loss function.   Let's starts from the origin VPFP equation \eqref{newvfp}:
\begin{equation*}
    \begin{aligned}
        \partial_t f+v \cdot \nabla_x f &=\frac{1}{\varepsilon} \nabla_v\cdot\left[(v+\nabla_x\phi) f+\nabla_v f\right]=: \frac{1}{\varepsilon} \mathscr{L} f,\\
        -\Delta_x \phi &=\rho-h(x).
    \end{aligned}
\end{equation*}
Due to the property of mass conservation of Fokker-Planck-type operator $\mathscr{L}$, one can integrate the first equation over $v$ to obtain the equation of mass conservation law:
\begin{equation}\label{rho evl}
    \partial_t \rho+\nabla_x\cdot\langle vf\rangle=0.
\end{equation}
Our new model is the Vlasov equation \eqref{vlasov} together with an extra condition: the mass conservation law
\begin{equation}\label{system2}
    \left\{\begin{aligned}
        \partial_t f+v\cdot \nabla_x f &= \frac{1}{\varepsilon} \mathscr{L} f,\\
        \partial_t \rho+\nabla_x\cdot\langle vf\rangle &=0,\\
        -\triangle_x\phi &=\rho-h,\\
        \rho &= \langle f\rangle.
    \end{aligned}\right.
\end{equation}
This system  is equivalent to the original system \eqref{vlasov}-\eqref{poisson}, since equation \eqref{rho evl} is just the integration of the Vlasov equation \eqref{newvfp} which is automatically satisfied in the continuous model but not necessarily so for DNN approximations. In the new model, the concerned macroscopic physical quantity $\rho$ is naturally incorporated by the conservation law without any subsidiary conditions. And we will see the high-field limit equation naturally follows from the new model as $\varepsilon \to 0$. 

The formal derivation of the limit equation starts from the new model \eqref{system2}. Multiplying the first equation of the model \eqref{system2} by $\varepsilon$ and yielding $\varepsilon \rightarrow 0$, this immediately gives $\mathscr{L} f=0$. Since $f\in\mathscr{N}(\mathscr{L})=\{f=\rho \mathscr{M}$, where $\rho:=\langle f\rangle\}$, it follows that $f=\rho \mathscr{M}$. Substituting this form of $f$ into the second equation of the model \eqref{system2}, and repeating the derivations in Section 1.3, limit system degenerates to:
\begin{equation}
    \left\{\begin{array}{l}
    \partial_t \rho-\nabla_x \cdot\left(\rho \nabla_x \phi\right)=0,\\
    -\triangle_x \phi=\rho-h(x).
    \end{array}\right.
\end{equation}

Furthermore, one takes care of the corresponding initial-boundary value problem for the new model \eqref{system2}. The initial and boundary conditions of $f$ are directly extended to $\rho$, since $\rho$ is direct integration of $f$ over the velocity space.

\subsubsection{A new loss function with AP property}

In this part, a new APNN method is proposed to solve the VPFP system. The key approach is incorporating the mass conservation based model \eqref{system2}, which possess the AP property, into the loss of the vanilla PINNs. The starting initial-boundary problem over a bounded domain is:
\begin{equation}
    \begin{cases}
        \partial_t f+v \cdot \nabla_x f=\frac{1}{\varepsilon} \mathscr{L} f, & (t, x, v) \in \mathcal{T} \times \mathcal{D} \times \Omega\\
        \partial_t \rho+\nabla_x\cdot\langle vf\rangle =0, & (t, x) \in \mathcal{T} \times \partial \mathcal{D}\\
        -\triangle_x \phi(t, x) =\rho(t, x)-h(x), & (t, x) \in \mathcal{T} \times \mathcal{D},\\
        \mathcal{B} f=F_{\mathrm{B}}, & (t, x, v) \in \mathcal{T} \times \partial \mathcal{D} \times \Omega\\ 
        \mathcal{B} \rho= P_{\mathrm{B}}, & (t, x) \in \mathcal{T} \times \partial \mathcal{D}\\
        \mathcal{B} \phi= \Phi_{\mathrm{B}}, & (t, x) \in \mathcal{T} \times \partial \mathcal{D}\\
        \mathcal{I} f=f_0, & (t, x, v) \in\{t=0\} \times \mathcal{D} \times \Omega\\
        \mathcal{I} \rho=\rho_0, & (t, x) \in\{t=0\} \times \mathcal{D}\\
        \mathcal{I} \phi=\phi_0, & (t, x) \in\{t=0\} \times \mathcal{D}
    \end{cases}
\end{equation}
where $F_{\mathrm{B}}, f_0, P_{\mathrm{B}}, \rho_0, \Phi_{\mathrm{B}}, \phi_0$ are given functions. 

First, we parametrize the three functions $\rho$, $f$ and $\phi$ using DNNs, denoted by $\eta$ for brevity. More precisely, we set:
\begin{equation}
    \begin{aligned}
        \rho_\eta^{\mathrm{NN}}(t, x):= \log(1 + \exp \left(\tilde{\rho}_\eta^{\mathrm{NN}}(t, x)\right)) \approx \rho(t, x),\\
        f_\eta^{\mathrm{NN}}(t, x):= \log(1 + \exp \left(\tilde{f}_\eta^{\mathrm{NN}}(t, x)\right)) \approx f(t, x),
    \end{aligned}
\end{equation}
which preserve the non-negative property of $\rho$ and $f$. The function $\phi$ is parameterized as:
\begin{equation}
    \phi_\eta^{\mathrm{NN}}(t, x)\approx \phi(t, x).
\end{equation}
Here $\tilde{\rho}_\theta^{\mathrm{NN}}$, $\tilde{f}_\theta^{\mathrm{NN}}$ and $\phi_\theta^{\mathrm{NN}}$ are all fully-connected neural networks. Next, we introduce the least square of the residual of the mass conservation based model \eqref{system2} as a new APNN loss:
\begin{equation}\label{loss2}
    \begin{aligned}
        \mathcal{R}_{\mathrm{APNN}}^{\varepsilon}=& \frac{\kappa_1^\text{macro}}{|\mathcal{T} \times \mathcal{D}|} \int_{\mathcal{T}}\int_{\mathcal{D}}\left|\partial_t \rho_\eta^{\mathrm{NN}} + \partial_x\langle vf_\eta^{\mathrm{NN}}\rangle \right|^2\mathrm{~d} x \mathrm{d} t\\
        &+\frac{\kappa_1^\text{kinetic}}{|\mathcal{T} \times \mathcal{D} \times \Omega|} \int_{\mathcal{T}} \int_{\mathcal{D}} \int_{\Omega}\left|\varepsilon\partial_tf_\eta^{\mathrm{NN}}+\varepsilon v\partial_xf_\eta^{\mathrm{NN}}-\mathscr Lf_\eta^{\mathrm{NN}} \right|^2\mathrm{~d}v\mathrm{~d} x \mathrm{d} t\\
        &+ \frac{\kappa_1^\text{Poisson}}{|\mathcal{T} \times \mathcal{D}|} \int_{\mathcal{T}}\int_{\mathcal{D}}\left|-\partial_x^2\phi_\eta^{\mathrm{NN}}-(\rho_\eta^{\mathrm{NN}}-h)\right|^2\mathrm{~d} x \mathrm{d} t \\
        &+\frac{\kappa_2^\rho}{|\mathcal{T} \times \partial \mathcal{D}|} \int_{\mathcal{T}} \int_{\partial \mathcal{D}} \left|\mathcal{B}\rho_\eta^{\mathrm{NN}}-P_{\mathrm{B}}\right|^2 \mathrm{d} x \mathrm{d} t \\
        &+\frac{\kappa_2^f}{|\mathcal{T} \times \partial \mathcal{D} \times \Omega\mid} \int_{\mathcal{T}} \int_{\partial \mathcal{D}} \int_{\Omega}\left|\mathcal{B}f_\eta^{\mathrm{NN}}-F_{\mathrm{B}}\right|^2 \mathrm{~d} v \mathrm{d} x \mathrm{d} t \\
        &+\frac{\kappa_2^\phi}{|\mathcal{T} \times \partial \mathcal{D}|} \int_{\mathcal{T}} \int_{\partial \mathcal{D}}\left|\mathcal{B}\phi_\eta^{\mathrm{NN}}-\Phi_{\mathrm{B}}\right|^2 \mathrm{d} x \mathrm{d} t \\
        &+\frac{\kappa_3^\rho}{|\mathcal{D}|} \int_{\mathcal{D}} \left|\mathcal{I}\rho_\eta^{\mathrm{NN}}-\rho_0\right|^2\mathrm{d} x
        +\frac{\kappa_3^f}{|\mathcal{D} \times \Omega|} \int_{\mathcal{D}} \int_{\Omega}\left|\mathcal{I}f_\eta^{\mathrm{NN}}-f_0\right|^2 \mathrm{~d} v \mathrm{d} x\\ 
        &+\frac{\kappa_3^\phi}{|\mathcal{D}|} \int_{\mathcal{D}} \left|\mathcal{I}\phi_\eta^{\mathrm{NN}}-\phi_0\right|^2\mathrm{d} x\\
        &+ \frac{\kappa_4}{|\mathcal{T} \times \mathcal{D}|} \int_{\mathcal{T}} \int_{\mathcal{D}}\left|\left\langle f_\eta^\mathrm{NN}\right\rangle-\rho_\eta^\mathrm{NN}\right|^2 \mathrm{~d} \boldsymbol{x} \mathrm{d} t
    \end{aligned}
\end{equation}
For ease of description, loss $\mathcal{R}_{\mathrm{APNN}}^{\varepsilon}$ is abbreviated as:
\begin{equation}
    \mathcal{R}_{\mathrm{APNN}}^{\varepsilon} := \kappa_1\mathcal{R}_{\text{residual}}^{\varepsilon} + \kappa_2\mathcal{R}_{\text{bc}}^{\varepsilon} + \kappa_3\mathcal{R}_{\text{ic}}^{\varepsilon}.
\end{equation}

At last, we show the AP property of the loss \eqref{loss2}. We just focus on the residual term in loss:
\begin{equation}
    \begin{aligned}
        \mathcal{R}_{\text{residual}}^{\varepsilon}=& \frac{1}{|\mathcal{T} \times \mathcal{D}|} \int_{\mathcal{T}}\int_{\mathcal{D}}\left|\partial_t \rho_\eta^{\mathrm{NN}} + \partial_x\langle vf_\eta^{\mathrm{NN}}\rangle \right|^2\mathrm{~d} x \mathrm{d} t\\
        &+\frac{1}{|\mathcal{T} \times \mathcal{D} \times \Omega|} \int_{\mathcal{T}} \int_{\mathcal{D}} \int_{\Omega}\left|\varepsilon\partial_tf_\eta^{\mathrm{NN}}+\varepsilon v\partial_xf_\eta^{\mathrm{NN}}-\mathscr Lf_\eta^{\mathrm{NN}} \right|^2\mathrm{~d}v\mathrm{~d} x \mathrm{d} t\\
        &+ \frac{1}{|\mathcal{T} \times \mathcal{D}|} \int_{\mathcal{T}}\int_{\mathcal{D}}\left|-\partial_x^2\phi_\eta^{\mathrm{NN}}-(\rho_\eta^{\mathrm{NN}}-h)\right|^2\mathrm{~d} x \mathrm{d} t .
    \end{aligned}
\end{equation}
Taking $\varepsilon \rightarrow 0$, formally it leads to:
\begin{equation}
    \begin{aligned}
        \mathcal{R}_{\text{residual}}^{\varepsilon}=& \frac{1}{|\mathcal{T} \times \mathcal{D}|} \int_{\mathcal{T}}\int_{\mathcal{D}}\left|\partial_t \rho_\eta^{\mathrm{NN}} + \partial_x\langle vf_\eta^{\mathrm{NN}}\rangle \right|^2\mathrm{~d} x \mathrm{d} t\\
        &+\frac{1}{|\mathcal{T} \times \mathcal{D} \times \Omega|} \int_{\mathcal{T}} \int_{\mathcal{D}} \int_{\Omega}\left|\mathscr Lf_\eta^{\mathrm{NN}} \right|^2\mathrm{~d}v\mathrm{~d} x \mathrm{d} t\\
        &+ \frac{1}{|\mathcal{T} \times \mathcal{D}|} \int_{\mathcal{T}}\int_{\mathcal{D}}\left|-\partial_x^2\phi_\eta^{\mathrm{NN}}-(\rho_\eta^{\mathrm{NN}}-h)\right|^2\mathrm{~d} x \mathrm{d} t .
    \end{aligned}
\end{equation}
This is the least square residual of the system:
\begin{equation}
    \left\{\begin{aligned}
    &\mathscr{L} f = 0, \\
    &\partial_t \rho+\partial_x\langle vf\rangle =0, \\
    &-\partial_x^2 \phi =\rho-h.
    \end{aligned}\right.
\end{equation}
The first equation $\mathscr{L}f = 0$ implies $f \in \mathscr{N}(\mathscr{L})$, inducing $f=\rho \mathscr{M}$. Plugging $f=\rho \mathscr{M}$ into the second equation and combining with the Poisson equation, one arrives at the limit system \eqref{limit1}. Therefore, the AP property of the new APNN loss has been concluded. 

Broadly, the mass conservation based APNN method holds several advantageous characteristics. First, the macroscopic physics quantity of interest is judiciously incorporated by the mass conservation law, free from restrictive assumptions. Secondly, the explicit Maxwellian form is rendered an unnecessary component.  A wider applicability is promised, since not all kinetic-type equations possess explicit Maxwellians. Moreover, Our new method is applicable to both equilibrium and non-equilibrium initial data.

\begin{remark}
For the mass conservation based APNN method, the mass conservation condition $\rho=\langle f\rangle$ is imposed to the loss as a soft constraint. This is different from the micro-macro decomposition based APNN method, in which the equivalent conservation condition $\langle g\rangle=0$ is built as a hard constraint in DNNs. As reported in \cite{jin2013asymptotic}, the constraint-type for the conservation mechanism significantly effects the training accuracy of the DNNs. We will examine these two processing methods for the conservation conditions in the numerical experiment section.
\end{remark}

%%%%%%%%%%%%%%%%%%%% Empirical loss %%%%%%%%%%%%%%%%%%%%%
\subsection{Empirical loss functions}

Since the loss functions for APNNs and PINNs are defined by integrals, in practice, one  needs to define the empirical loss as a Monte Carlo approximation to the integral loss. This involves the random selection to a small number of sub-domains (known as batches) to estimate the high-dimensional integrals. For the velocity-dependent operators $\langle\cdot\rangle$ and $\Pi(\cdot)$, the Gauss-Legendre quadrature provides additional accuracy. In detail, given $\left\{w_i, v_i^{\prime}\right\}_{i=1}^n$ are the nodes and weights, which $\left\{v_i^{\prime}\right\}_{i=1}^n$ are the roots of the Legendre polynomials of degree $n$ and $\left\{w_i\right\}_{i=1}^n$ are the corresponding weights. Then an integral can be approximated, for example, $\int_{-1}^1 f(v) \mathrm{d} v$ by the summation of linear combination of $f\left(v_i^{\prime}\right): \sum_{i-1}^n w_i f\left(v_i^{\prime}\right)$. 

We start with   the formulation of the empirical losses for vanilla PINN \eqref{loss PINN}  defined as:
\begin{equation}
    \begin{aligned}
        \mathcal{R}_{\mathrm{APNN}}^{\varepsilon}= &\frac{\mu_1^\text{Vlasov}}{L_1}\sum_{l=1}^{L_1}\left|\varepsilon\partial_tf_\gamma^{\mathrm{NN}}(t_l^{r}, x_l^{r}, v_l^{r})+\varepsilon v_l^{r}\partial_xf_\gamma^{\mathrm{NN}}(t_l^{r}, x_l^{r}, v_l^{r})-\mathscr Lf_\gamma^{\mathrm{NN}}(t_l^{r}, x_l^{r}, v_l^{r})\right|^2\\
        &+ \frac{\mu_1^\text{Poisson}}{L_1}\sum_{l=1}^{L_1}\left|-\partial_x^2\phi_\gamma^{\mathrm{NN}}(t_l^{r}, x_l^{r})-(\rho_\gamma^{\mathrm{NN}}(t_l^{r}, x_l^{r})-h(x_l^{r}))\right|^2\\
        &+\frac{\mu_2^f}{L_2}\sum_{l=1}^{L_2}\left|\mathcal{B}f_\gamma^{\mathrm{NN}}(t_l^{b}, x_l^{b}, v_l^{b})-F_{\mathrm{B}}(t_l^{b}, x_l^{b}, v_l^{b})\right|^2 + \frac{\mu_2^\phi}{L_2}\sum_{l=1}^{L_2}\left|\mathcal{B}\phi_\gamma^{\mathrm{NN}}(t_l^{b}, x_l^{b})-\Phi_{\mathrm{B}}(t_l^{b}, x_l^{b})\right|^2\\
        & + \frac{\mu_3^f}{L_3}\sum_{l=1}^{L_3}\left|\mathcal{I}f_\gamma^{\mathrm{NN}}(t_l^{k}, x_l^{k}, v_l^{k})-f_0(t_l^{k}, x_l^{k}, v_l^{k})\right|^2 + \frac{\mu_3^\phi}{L_3}\sum_{l=1}^{L_3}\left|\mathcal{I}\phi_\gamma^{\mathrm{NN}}(t_l^{k}, x_l^{k})-\phi_0(t_l^{k}, x_l^{k})\right|^2\\
    \end{aligned}
\end{equation}
where $\{(t_l^{r}, x_l^{r}, v_l^{r})\}_{l=1}^{L_1}\in\mathcal{T} \times \mathcal{D} \times \Omega$ is a batch of sample data for the residual, $\{(t_l^{b}, x_l^{b}, v_l^{b})\}_{l=1}^{L_2}\in\mathcal{T} \times \partial \mathcal{D} \times \Omega$ is a batch of sample data, and $\{(t_l^{k}, x_l^{k}, v_l^{k})\}_{l=1}^{L_3}\in\{0\} \times \mathcal{D} \times \Omega$ is a batch of sample data for initial conditions. 

Similarly, the empirical loss function for the micro-macro decomposition based APNN method \eqref{MMloss} is defined as:
\begin{equation}
    \begin{aligned}
    \mathcal{R}_{\mathrm{APNN}}^{\varepsilon}= & \frac{\lambda_1^\text{macro}}{N_1} \sum_{i=1}^{N_1}\left|\partial_t \rho_\theta^{\mathrm{NN}}(t_i^{r}, x_i^{r}) + \partial_x\langle v(\rho_\theta^{\mathrm{NN}}\mathscr M_\theta^{\mathrm{NN}})\rangle(t_i^{r}, x_i^{r}) +\varepsilon\partial_x\langle vg_\theta^{\mathrm{NN}}\rangle (t_i^{r}, x_i^{r})\right|^2 \\
    & +\frac{\lambda_1^\text{micro}}{N_1} \sum_{i=1}^{N_1} \left| \varepsilon\partial_t g_\theta^{\mathrm{NN}}(t_i^{r}, x_i^{r}, v_i^{r}) + \varepsilon(v_i^{r}\partial_x g_\theta^{\mathrm{NN}}(t_i^{r}, x_i^{r}, v_i^{r})-\partial_x\langle vg_\theta^{\mathrm{NN}}\rangle(t_i^{r}, x_i^{r})\mathscr M_\theta^{\mathrm{NN}}(t_i^{r}, x_i^{r}, v_i^{r})) \right. \\
    &\left. - \{\mathscr{L} g_\theta^{\mathrm{NN}}(t_i^{r}, x_i^{r}, v_i^{r})-v_i^{r} \partial_x(\rho_\theta^{\mathrm{NN}}(t_i^{r}, x_i^{r}) \mathscr{M}_\theta^{\mathrm{NN}}(t_i^{r}, x_i^{r}, v_i^{r}))\right.\\
    & \left.+\partial_x\langle v(\rho_\theta^{\mathrm{NN}} \mathscr{M}_\theta^{\mathrm{NN}})\rangle (t_i^{r}, x_i^{r}) \mathscr{M}_\theta^{\mathrm{NN}}(t_i^{r}, x_i^{r}, v_i^{r})-\rho_\theta^{\mathrm{NN}}(t_i^{r}, x_i^{r})(\partial_t\mathscr{M}_\theta^{\mathrm{NN}}(t_i^{r}, x_i^{r}, v_i^{r}))\}\right|^2 \\
    & \frac{\lambda_1^\text{Poisson}}{N_1} \sum_{i=1}^{N_1}\left|-\partial_x^2\phi_\theta^{\mathrm{NN}}(t_i^{r}, x_i^{r})-(\rho_\theta^{\mathrm{NN}}(t_i^{r}, x_i^{r})-h(x_i^{r}))\right|^2 \\
    & +\frac{\lambda_2^f}{N_2} \sum_{i=1}^{N_2}\left|\mathcal{B}\left(\rho_\theta^{\mathrm{NN}}(t_i^b, x_i^b)+\varepsilon g_\theta^{\mathrm{NN}}(t_i^b, x_i^b, v_i^b)\right)-F_{\mathrm{B}}(t_i^b, x_i^b, v_i^b)\right|^2 \\
    & +\frac{\lambda_2^\phi}{N_2} \sum_{i=1}^{N_2}\left|\mathcal{B}\phi_\theta^{\mathrm{NN}}(t_i^b, x_i^b)-\Phi_{\mathrm{B}}(t_i^b, x_i^b)\right|^2 \\
    & +\frac{\lambda_3^f}{N_3} \sum_{i=1}^{N_3}\left|\mathcal{I}\left(\rho_\theta^{\mathrm{NN}}\left(t_i^k, x_i^k\right)+\varepsilon g_\theta^{\mathrm{NN}}\left(t_i^k, x_i^k, v_i^k\right)\right)-f_0\left(t_i^k, x_i^k, v_i^k\right)\right|^2\\
    & +\frac{\lambda_3^\phi}{N_3} \sum_{i=1}^{N_3}\left|\mathcal{I}\phi_\theta^{\mathrm{NN}}\left(t_i^k, x_i^k\right)-\phi_0\left(t_i^k, x_i^k\right)\right|^2
    \end{aligned}
\end{equation}
where $\{(t_i^{r}, x_i^{r}, v_i^{r})\}_{i=1}^{N_1}\in\mathcal{T} \times \mathcal{D} \times \Omega$ is a batch of sample data for the residual, $\{(t_i^{b}, x_i^{b}, v_i^{b})\}_{i=1}^{N_2}\in\mathcal{T} \times \partial \mathcal{D} \times \Omega$ is a batch of sample data for boundary conditions, and $\{(t_i^{k}, x_i^{k}, v_i^{k})\}_{i=1}^{N_3}\in\{0\} \times \mathcal{D} \times \Omega$ is a batch of sample data for initial conditions.

The empirical loss function for the mass conservation based APNN method \eqref{loss2} is analogously defined as:
\begin{equation}
    \begin{aligned}
        \mathcal{R}_{\mathrm{APNN}}^{\varepsilon}=& \frac{\kappa_1^\text{macro}}{M_1}\sum_{j=1}^{M_1}\left|\partial_t \rho_\eta^{\mathrm{NN}}(t_j^{r}, x_j^{r}) + \partial_x\langle vf_\eta^{\mathrm{NN}}(t_j^{r}, x_j^{r})\rangle\right|^2\\
        &+ \frac{\kappa_1^\text{kinetic}}{M_1}\sum_{j=1}^{M_1}\left|\varepsilon\partial_tf_\eta^{\mathrm{NN}}(t_j^{r}, x_j^{r}, v_j^{r})+\varepsilon v_j^{r}\partial_xf_\eta^{\mathrm{NN}}(t_j^{r}, x_j^{r}, v_j^{r})-\mathscr Lf_\eta^{\mathrm{NN}}(t_j^{r}, x_j^{r}, v_j^{r})\right|^2\\
        &+ \frac{\kappa_1^\text{Poisson}}{M_1}\sum_{j=1}^{M_1}\left|-\partial_x^2\phi_\eta^{\mathrm{NN}}(t_j^{r}, x_j^{r})-(\rho_\eta^{\mathrm{NN}}(t_j^{r}, x_j^{r})-h(x_j^{r}))\right|^2\\
        &+\frac{\kappa_2^\rho}{M_2}\sum_{j=1}^{M_2}\left|\mathcal{B}\rho_\eta^{\mathrm{NN}}(t_j^{b}, x_j^{b})-P_{\mathrm{B}}(t_j^{b}, x_j^{b})\right|^2+\frac{\kappa_2^f}{M_2}\sum_{j=1}^{M_2}\left|\mathcal{B}f_\eta^{\mathrm{NN}}(t_j^{b}, x_j^{b}, v_j^{b})-F_{\mathrm{B}}(t_j^{b}, x_j^{b}, v_j^{b})\right|^2\\
        &+ \frac{\kappa_2^\phi}{M_2}\sum_{j=1}^{M_2}\left|\mathcal{B}\phi_\eta^{\mathrm{NN}}(t_j^{b}, x_j^{b})-\Phi_{\mathrm{B}}(t_j^{b}, x_j^{b})\right|^2\\
        &+\frac{\kappa_3^\rho}{M_3}\sum_{j=1}^{M_3}\left|\mathcal{I}\rho_\eta^{\mathrm{NN}}(t_j^{k}, x_j^{k})-\rho_0(t_j^{k}, x_j^{k})\right|^2 + \frac{\kappa_3^f}{M_3}\sum_{j=1}^{M_3}\left|\mathcal{I}f_\eta^{\mathrm{NN}}(t_j^{k}, x_j^{k}, v_j^{k})-f_0(t_j^{k}, x_j^{k}, v_j^{k})\right|^2\\
        &+ \frac{\kappa_3^\phi}{M_3}\sum_{j=1}^{M_3}\left|\mathcal{I}\phi_\eta^{\mathrm{NN}}(t_j^{k}, x_j^{k})-\phi_0(t_j^{k}, x_j^{k})\right|^2\\
        &+ \frac{\kappa_4}{M_4}\sum_{j=1}^{M_4}\left|\left\langle f_\eta^\mathrm{NN}\right\rangle(t_j^{c}, x_j^{c})-\rho_\eta^\mathrm{NN}(t_j^{c}, x_j^{c})\right|^2 
    \end{aligned}
\end{equation}
where $\{(t_j^{r}, x_j^{r}, v_j^{r})\}_{j=1}^{M_1}\in\mathcal{T} \times \mathcal{D} \times \Omega$ is a batch of sample data for the residual, $\{(t_j^{b}, x_j^{b}, v_j^{b})\}_{j=1}^{M_2}\in\mathcal{T} \times \partial \mathcal{D} \times \Omega$ is  a batch of sample data for boundary conditions, $\{(t_j^{k}, x_j^{k}, v_j^{k})\}_{j=1}^{M_3}\in\{0\} \times \mathcal{D} \times \Omega$ is a batch of sample data for initial conditions, and $\{(t_j^{c}, x_j^{c})\}_{j=1}^{M_4}\in\mathcal{T} \times \mathcal{D}$ is a batch of sample data for the mass conservation.

 With the empirical losses formed, one can generate the training set, and apply the optimization procedure for the empirical losses to approximate the solution of the VPFP system.

%%%%%%%%%%%%%%%%%%%%%%%%%%%%%%%%%%%%%%%%%%%%%%%%
%%%%%%%%%%%%%%%%%%%%%%%%%%%%%%%%%%%%%%%%%%%%%%%%
\section{Numerical examples}

In this section, extensive numerical results are presented for several problems chosen from kinetic regimes $(\varepsilon \approx O(1))$ to high-field regimes $(\varepsilon \rightarrow 0)$, in order to verify two proposed APNN methods and compare their performances.

Some settings for the numerical experiments are stated as follows. The Adam version of the gradient descent algorithm is used for optimization of the loss functions. Neural network parameters $\theta$, $\eta$ and $\gamma$ are initialized by the Xavier initialization. Fully-connected networks with $5$ hidden layers and $\tanh$ activation function are used. The number of Gauss-quadrature points is set to be $32$ for integration calculation. All hyper-parameters have been carefully tuned for optimal performance.

The reference solutions are obtained by standard finite difference methods. We will check the relative $\ell^2$ errors of several macroscopic quantities, such as the density $\rho(t, x)$ and the electric field $E(t, x) = -\partial_x\phi(t, x)$, between the results of DNN methods and reference solutions. The root mean square errors (RMSE for brevity) are also provided when necessary. The specific definition of these two errors of $\rho$ as an example are:
\begin{equation}
    \text{RMSE}(\rho):=\sqrt{\frac{1}{N}\sum_{j=1}^N\left|\rho_{\theta, j}^{\mathrm{NN}}-\rho_j^{\mathrm{ref}}\right|^2}, \quad
    \ell^2(\rho):=\sqrt{\frac{\sum_j\left|\rho_{\theta, j}^{\mathrm{NN}}-\rho_j^{\mathrm{ref}}\right|^2}{\sum_j\left|\rho_j^{\mathrm{ref}}\right|^2}},
\end{equation}
where $N$ is the total number of mesh points.

%%%%%%%%%%%%%%%%%%%% LD %%%%%%%%%%%%%%%%%%%%%
\subsection{Problem I: Landau damping}

In this part, the Landau damping case is exhibited to verify the AP property of the two proposed APNN loss functions from kinetic regimes to high-field regimes and compare their performances. The failure of vanilla PINN method in high-field regime is also demonstrated. The initial condition near the equilibrium takes the form:
\begin{equation}\label{Landau1}
    f_0(x, v)=\frac{1}{\sqrt{2 \pi}}\rho_0(x) \exp \left(-\frac{v^2}{2}\right), \quad(x, v) \in[0,2 \pi / k] \times \mathbb{R},
\end{equation}
where the initial density $\rho_0$ is:
\begin{equation}\label{Landau2}
    \rho_0(x) = 1+\alpha \cos (k x).
\end{equation}
The wave number $k$ is set to be $0.5$ and the amplitude of the perturbation $\alpha$ is set as $0.05$. The computational domain is $[0,2 \pi / k] \times\left[v_{\min }, v_{\max }\right],-v_{\min }=v_{\max }=6$. The initial electric field $\phi^0$ is the solution to the Poisson equation with $h(x)=1$. It is solved analytically as:
\begin{equation}
    \phi_0(x) = \frac{\alpha}{k^2}\cos(kx).
\end{equation}
Periodic boundary conditions are applied in the $x$ domain.

The detailed initial conditions (ICs) and boundary conditions (BCs) are specified for the APNN approaches. For the micro-macro decomposition based APNN method, the ICs for $\rho_0, g_0, \phi_0$ are given by:
\begin{equation}
    \rho_0(x) = 1+\alpha\cos(kx), \quad \phi_0(x) = \frac{\alpha}{k^2}\cos(kx), \quad g_0(x,v)=\frac{1}{\varepsilon}\left\{f_0(x,v)-\rho_0(x) \mathscr{M}_0(x,v)\right\},
\end{equation}
where $\mathscr{M}_0$ is the Maxwellian as: 
\begin{equation}
    \mathscr{M}_0(x, v)=\frac{1}{\sqrt{2\pi}}\exp\left(-\frac{(v+\partial_x\phi_0(x))^2}{2}\right).
\end{equation}
Periodic BCs on $\rho, g$ are directly imposed based on the periodic BC for $f$. For mass conservation based APNN method and vanilla PINN method, the IC and periodic BC on $f$ are immediately translated into the related settings. To improve the numerical performance, exact periodic BCs are enforced. The ansatz is based on a Fourier basis. We also construct a transform $\mathscr{T}: x\rightarrow\{\sin (kjx), \cos (kjx)\}_{j=1}^n$ before the first layer of neural networks and set $n=1$.

\begin{figure}[htbp]
	\centering  %图片全局居中
	\subfigbottomskip=4pt %两行子图之间的行间距
	\subfigcapskip=-8pt %设置子图与子标题之间的距离
	\subfigure[Kinetic regime $(\varepsilon=1)$.]{
		\includegraphics[width=0.8\linewidth]{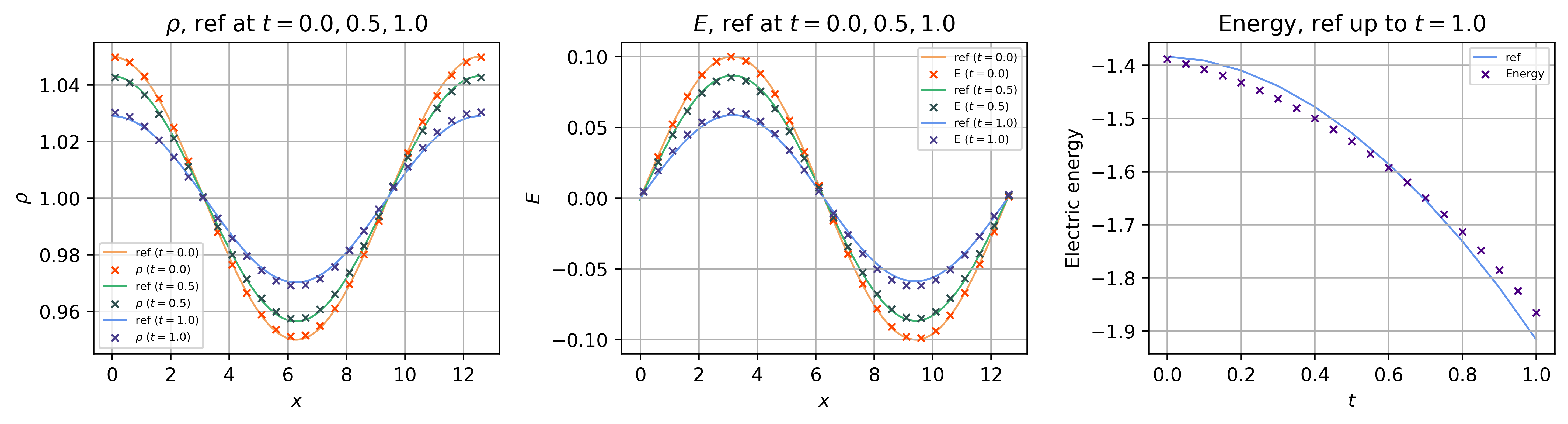}}
	\subfigure[Intermediate regime $(\varepsilon=0.5)$.]{
		\includegraphics[width=0.8\linewidth]{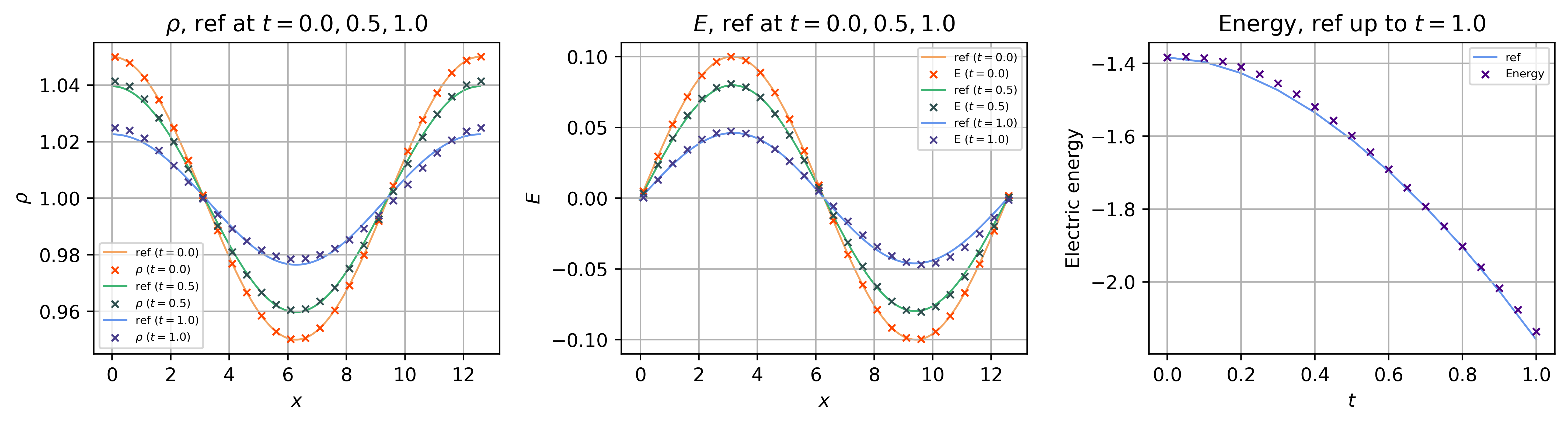}}
	\subfigure[High-field regime $(\varepsilon=0.01)$.]{
		\includegraphics[width=0.8\linewidth]{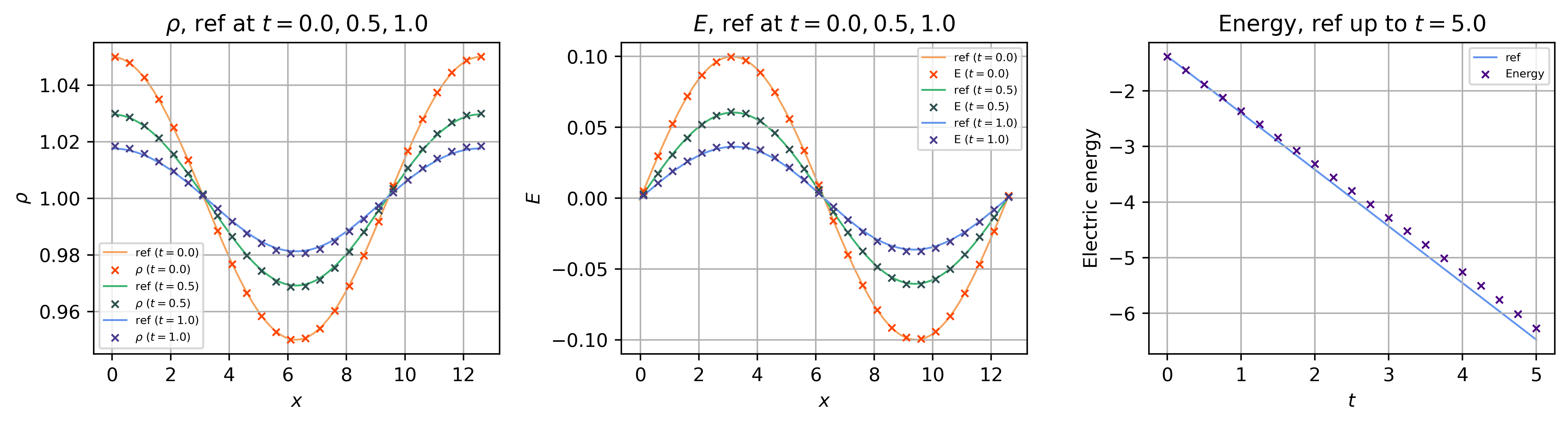}}
	\caption{Multiscale linear Landau damping solved by the micro-macro decomposition based APNN method. Density $\rho$ (left column) and electric field $E$ (middle column) as functions of space $x$ at $t = 0.0, 0.5, 1.0$; electric energy (right column) as a function of time $t$ up to $t=1.0$ with the kinetic and intermediate regimes, and $t=5.0$ with the high-field regime. The electric energy is plotted in log scale. Neural networks are $[3, 128, 128, 128, 128, 128, 1]$ for $\rho, \phi$ and $[4, 256, 256, 256, 256, 256, 1]$ for $g$. Batch size is $512$ in domain and $256$ for initial condition. Penalty $\lambda_1 = 300$, $\lambda_3 = 1$ for kinetic regime; penalty $\lambda_1 = 0.5$, $\lambda_3 = 1$ for intermediate regime; and penalty $\lambda_1 = 0.5$, $\lambda_3 = 1$ for high-field regime.}
    \label{ld mm}
\end{figure}

\begin{figure}[htbp]
	\centering  %图片全局居中
	\subfigbottomskip=4pt %两行子图之间的行间距
	\subfigcapskip=-8pt %设置子图与子标题之间的距离
	\subfigure[Kinetic regime $(\varepsilon=1)$.]{
		\includegraphics[width=0.8\linewidth]{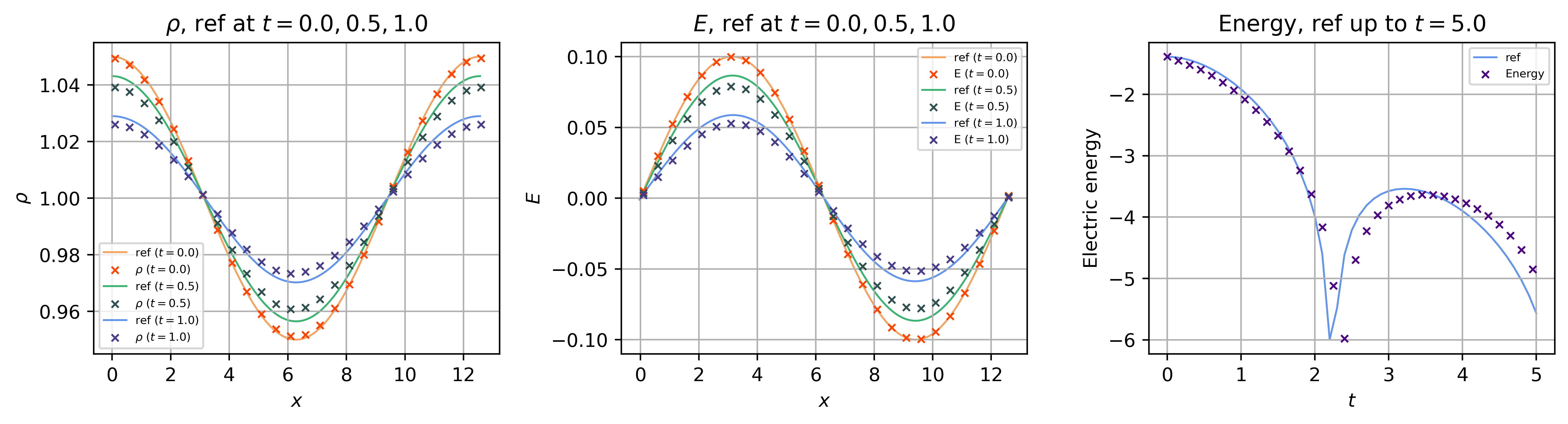}}
	\subfigure[Intermediate regime $(\varepsilon=0.5)$.]{
		\includegraphics[width=0.8\linewidth]{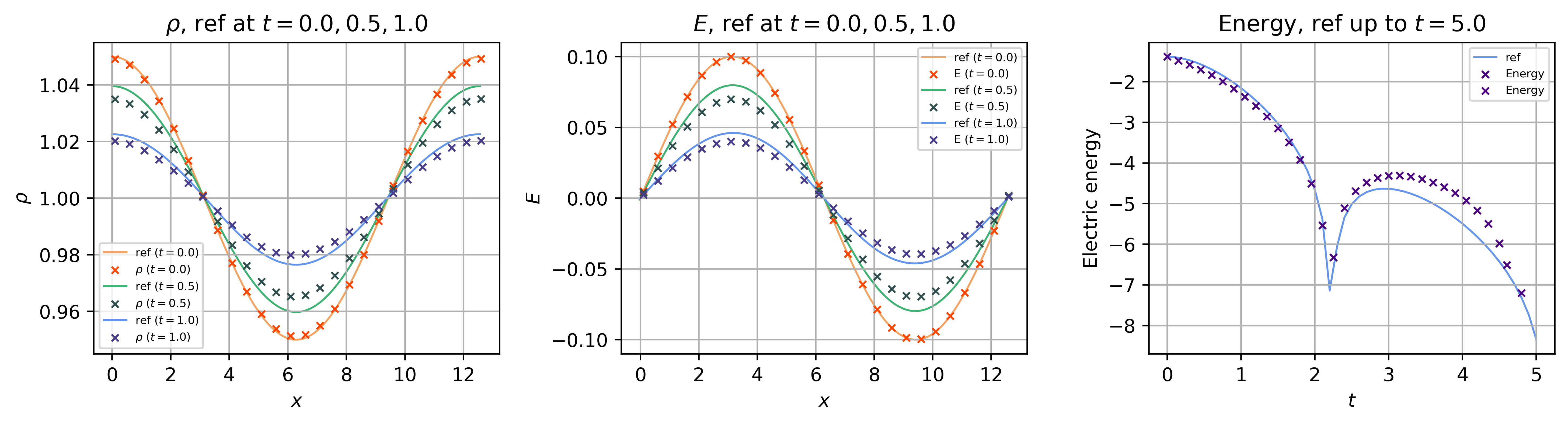}}
	\subfigure[High-field regime $(\varepsilon=0.01)$.]{
		\includegraphics[width=0.8\linewidth]{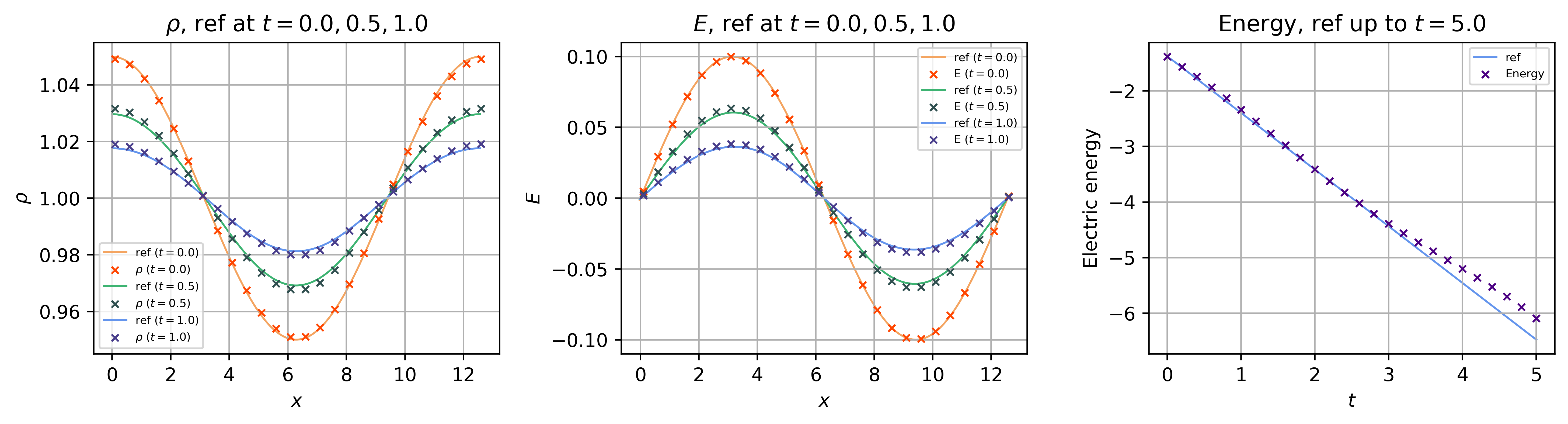}}
	\caption{Multiscale linear Landau damping solved by the mass conservation based APNN method. Density $\rho$ (left column) and electric field $E$ (middle column) as functions of space $x$ at $t = 0.0, 0.5, 1.0$; electric energy (right column) as a function of time $t$ up to $t=5.0$. Neural networks are $[3, 128, 128, 128, 128, 128, 1]$ for $\rho, \phi$ and $[4, 256, 256, 256, 256, 256, 1]$ for $f$. Batch size is $512$ in domain and $256$ for initial condition and $256$ for conservation condition $\rho=\langle f\rangle$. Penalty $\kappa_1 = 120$, $\kappa_3 = \kappa_4 = 1$ for kinetic regime; penalty $\kappa_1 = 150$, $\kappa_3 = \kappa_4 = 1$ for intermediate regime; and penalty $\kappa_1 = 500$, $\kappa_3 = \kappa_4 = 1$ for high-field regime.	}
    \label{ld mc}
\end{figure}

\begin{figure}[htbp]
	\centering  %图片全局居中
	\subfigbottomskip=4pt %两行子图之间的行间距
	\subfigcapskip=-8pt %设置子图与子标题之间的距离
	\subfigure[Kinetic regime $(\varepsilon=1)$.]{
		\includegraphics[width=0.8\linewidth]{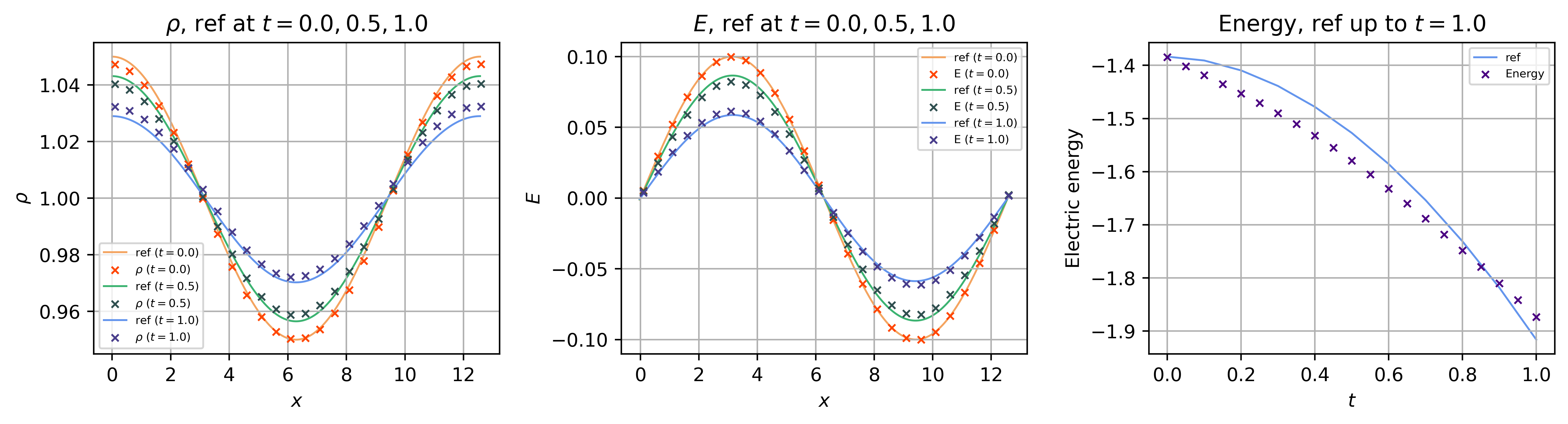}}
	\subfigure[Intermediate regime $(\varepsilon=0.5)$.]{
		\includegraphics[width=0.8\linewidth]{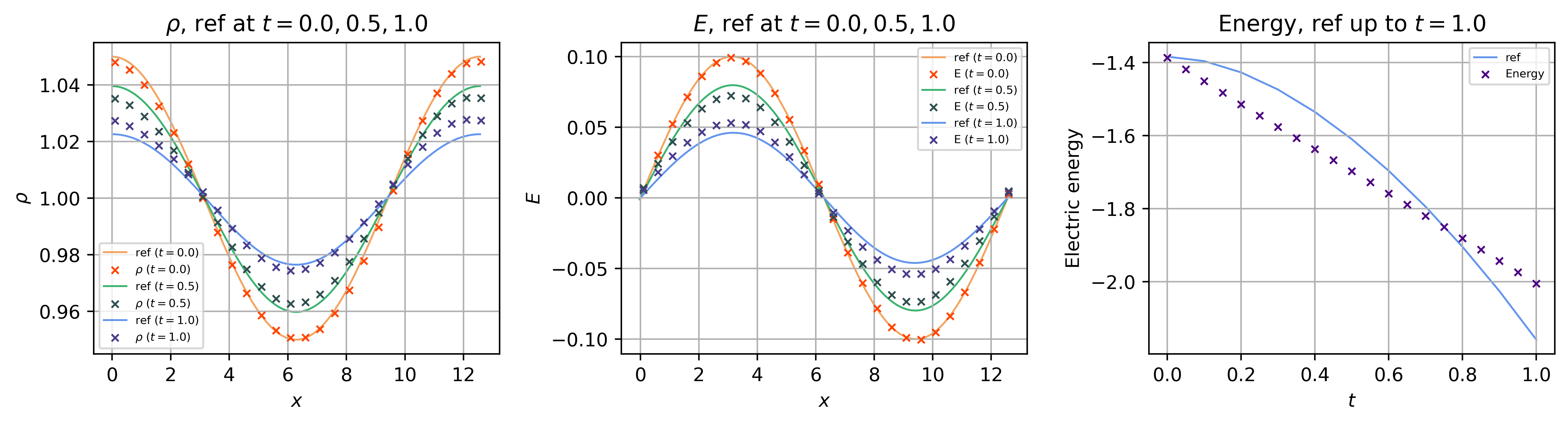}}
	\caption{Multiscale linear Landau damping solved by vanilla PINN method with kinetic and intermediate regimes. Density $\rho$ (column left) and electric field $E$ (column middle) as functions of space $x$ for $t = 0.0, 0.5, 1.0$; electric energy (right column) as a function of time $t$ up to $t=1.0$. Neural networks are $[3, 128, 128, 128, 128, 128, 1]$ for $\phi$ and $[4, 256, 256, 256, 256, 256, 1]$ for $f$. Batch size is $512$ in domain and $256$ for initial condition. Penalty $\mu_1 = 30$, $\mu_3 = 1$ for kinetic regime; and penalty $\mu_1 = 50$, $\mu_3 = 1$ for intermediate regime.}
    \label{ld pinn}
\end{figure}

\begin{figure}[htbp]
	\centering  %图片全局居中
	\subfigbottomskip=4pt %两行子图之间的行间距
	\subfigcapskip=-8pt %设置子图与子标题之间的距离
	\subfigure[Density $\rho$ (top) and electric field $E$ (bottom) as functions of space $x$ for $t = 0.0, 0.5, 1.0$.]{
		\includegraphics[width=0.8\linewidth]{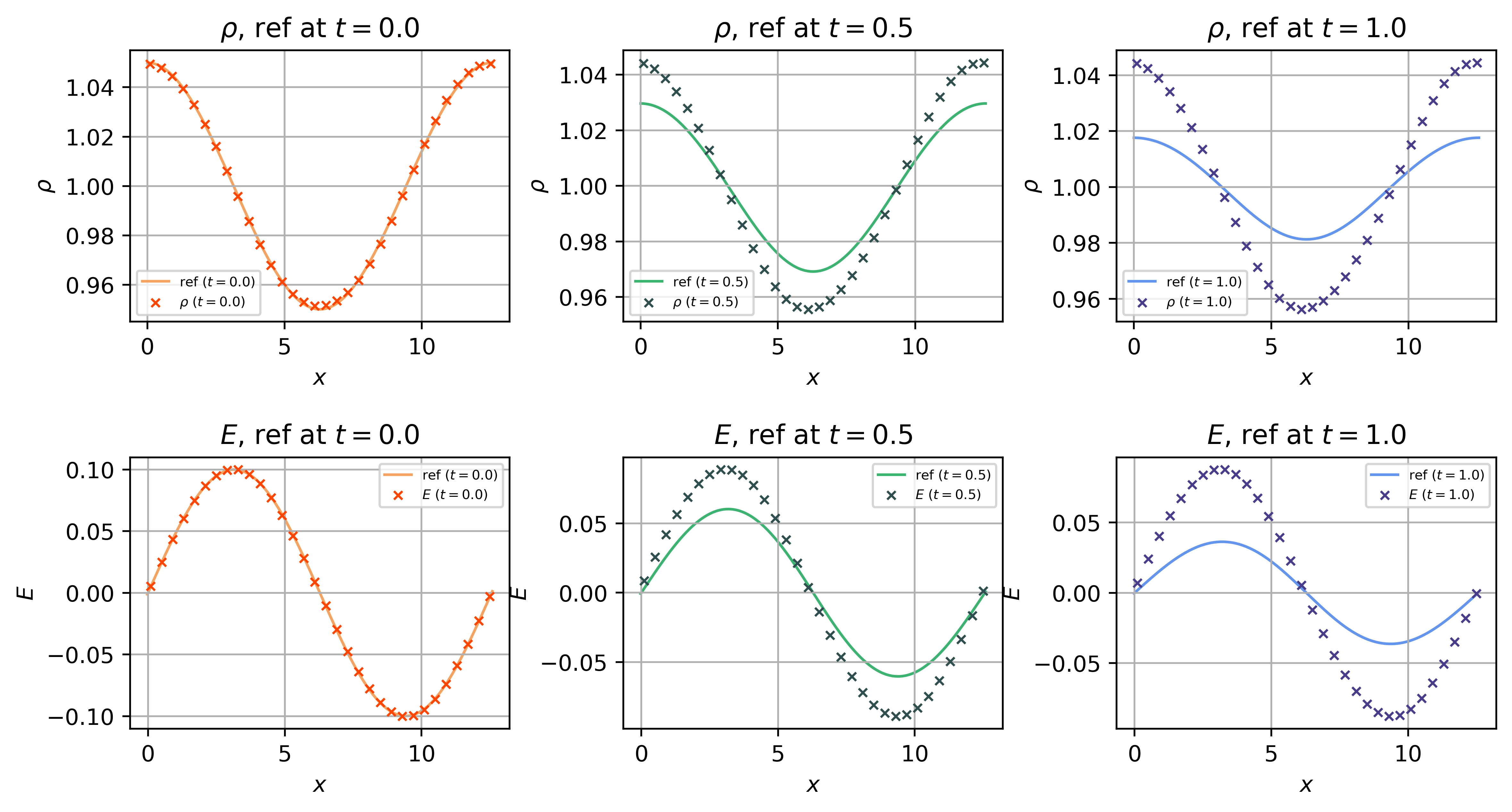}}
	\subfigure[Electric energy as a function of time $t$ up to $t=1.0$.]{
		\includegraphics[width=0.4\linewidth]{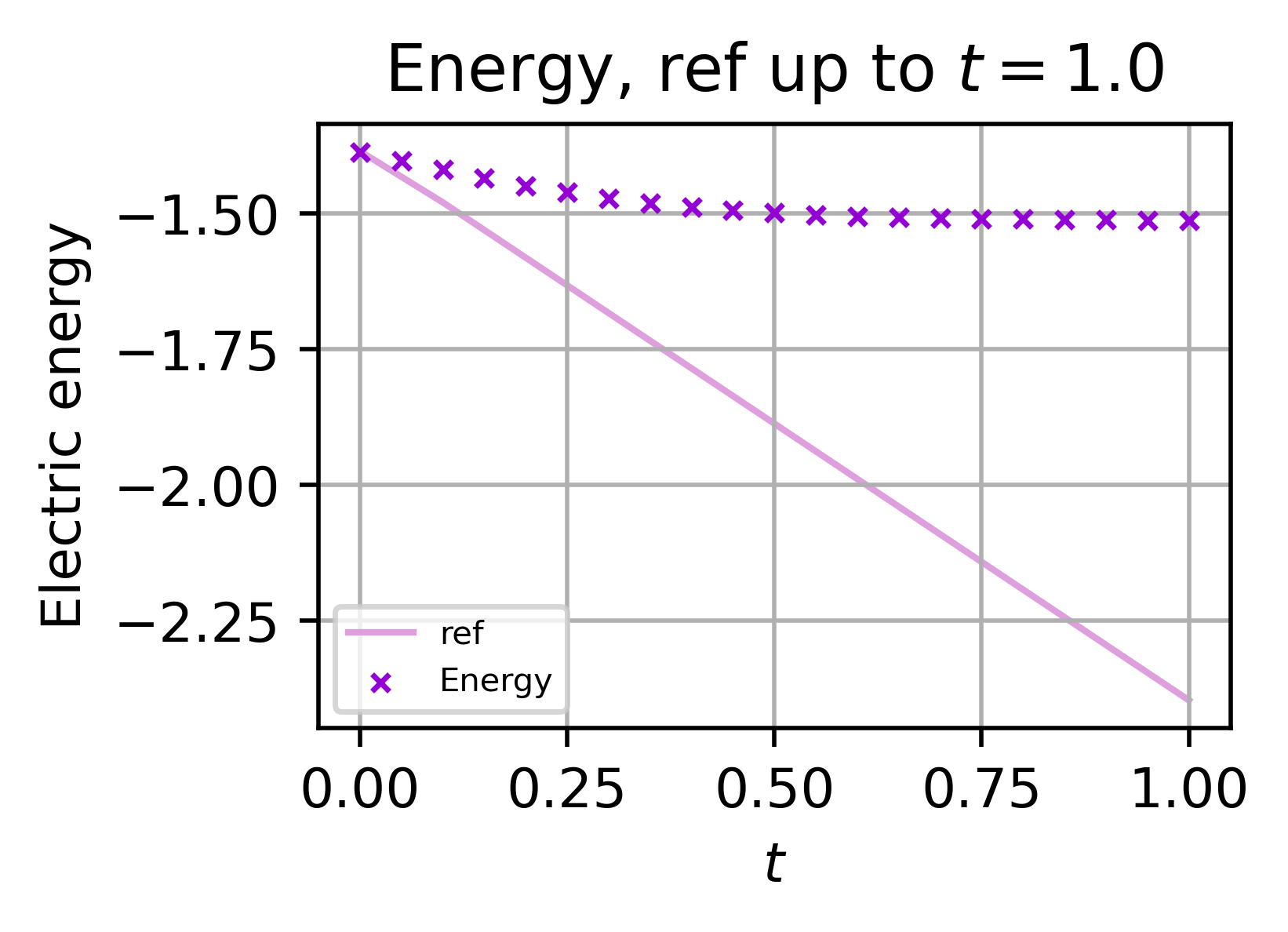}}
	\caption{Linear Landau damping solved by vanilla PINN method with high-field regime $(\varepsilon=0.01)$. Neural networks and batch size are the same as kinetic regime case. Penalty $\mu_1 = 100$, $\mu_3 = 1$.}
    \label{ld pinn high}
\end{figure}

Here, we present the performance of two APNN models in multiscale regimes, which are characterized by the parameter $\varepsilon$. For each APNN model, we plot the density $\rho$ and the electric field $E$ as a function of $x$ at different time, as well as the time evolution of the electric energy $\|E(t)\|_{L^2}$. Figure \ref{ld mm} and Figure \ref{ld mc} illustrate that the two proposed APNN models achieve good agreement with the reference solutions across kinetic, intermediate, and high-field regimes. These numerical results confirm the AP property of two proposed APNN methods and supports our analysis. As a comparison, we plot the result of vanilla PINN method in Figure \ref{ld pinn} and Figure \ref{ld pinn high}. The vanilla PINN method exhibits slightly poorer performance in kinetic regimes $(\varepsilon=1.0, 0.5)$ and fails to capture the asymptotic limit regime. Overall, two proposed APNN models improve the vanilla PINN method and have the capability to reliably reproduce solutions for multiscale plasma physics systems. 

Moreover, in the kinetic regime, we observe different performance of two proposed APNN methods at different time scales. The mass conservation based APNN method can capture the correct electric energy evolution over long duration up to $t=5.0$ (see right column in Figure \ref{ld mc}). However, the micro-macro decomposition based APNN method can only accurately characterize the dynamic behaviors in short duration up to $t=1.0$ (see errors reported in Table \ref{ld table kinetic}, Table \ref{ld table inter}, Table \ref{ld table high}). A possible reason is that the variable order of the non-equilibrium part $g$ may influence the training of DNNs \cite{lou2021physics}. As time increases, the magnitude of $g$ is significantly reduced. This may encumber the DNN training, thereby inducing the inaccuracy of the micro-macro decomposition based APNN method over long time duration. This weakness could be improved by employing self-adaptive PINNs \cite{mcclenny2020self, subramanian2022adaptive}. One promising strategy is to apply a time-dependent weight to non-equilibrium part $g$\cite{mcclenny2020self}. Another promising strategy is to allocate the time-dependent number of collocation points for non-equilibrium part $g$  \cite{subramanian2022adaptive}. These self-adaptive strategies could assist in precisely characterizing $g$ with variable order, thus improve the long time accuracy of the micro-macro decomposition based APNN method in VPFP system. In contrast, the mass conservation based APNN method gets rid of $g$ and introduces the density $\rho$. Density $\rho$ and distribution $f$ persist on equivalent first-orders, which broads its applicability over long time.

For the VPFP system, it is noteworthy that the hard or soft constraint-type of the mass conservation conditions do not significantly influence the accuracy of these two proposed APNN methods. In every numerical case tabled, the conservation condition ($\langle g\rangle=0$) is treated as a hard condition in the micro-macro decomposition based APNN method. We also examine the soft constraint ($\langle f\rangle=\rho$), constructed in the mass conservation based APNN method,  which has a performance similar to its hard constraint counterpart (see Figure \ref{ld mm} and Figure \ref{ld mc}). Thus, the APNN methods can achieve comparable performances with different constraint-type.

%%%%%%%%%%%%%% ld t1 %%%%%%%%%%%%%%%%%%%%%

\begin{minipage}[htbp]{\textwidth}
	\centering
	%subTable1
\begin{minipage}[!h]{\textwidth}
	\centering
% 	\makeatletter\def\@captype{}
	\begin{tabular}{ccc}
		\toprule
		$\ell^2(\rho)$ & $t=0.5$ & $t=1.0$ \\
		\midrule
		MM & $4.61\times 10^{-4}$  & $1.04\times 10^{-3}$ \\
		MC & $2.88\times 10^{-3}$ & $2.25\times 10^{-3}$ \\
		PINN & $1.60\times 10^{-3}$ & $2.48\times 10^{-3}$ \\
		\bottomrule
	\end{tabular}
	\quad
	\begin{tabular}{ccc}
		\toprule
		$\text{RMSE}(\rho)$ & $t=0.5$ & $t=1.0$ \\
		\midrule
		MM & $4.62\times 10^{-4}$  & $1.04\times 10^{-3}$ \\
		MC & $2.89\times 10^{-3}$ & $2.25\times 10^{-3}$  \\
		PINN & $1.61\times 10^{-3}$ & $2.48\times 10^{-3}$ \\
		\bottomrule
	\end{tabular}\\
	\vspace{0.02\textwidth}
    \begin{minipage}[!h]{\textwidth}
        \centering
        $(a)$ Errors of $\rho$ at $t=0.5, 1.0$.
    \end{minipage}
    \vspace{-0.005\textwidth}
\end{minipage}
    %
    %
    %subtable2
\begin{minipage}[htbp]{\textwidth}
	\centering
		\begin{tabular}{ccc}
		\toprule
		$\ell^2(E)$ & $t=0.5$ & $t=1.0$ \\
		\midrule
		MM & $1.64\times 10^{-2}$  & $5.71\times 10^{-2}$ \\
		MC & $9.36\times 10^{-2}$ & $1.10\times 10^{-1}$ \\
		PINN & $5.15\times 10^{-2}$ & $4.55\times 10^{-2}$ \\
		\bottomrule
	\end{tabular}
	\quad
	\begin{tabular}{ccc}
		\toprule
		$\text{RMSE}(E)$ & $t=0.5$ & $t=1.0$ \\
		\midrule
		MM & $1.01\times 10^{-3}$  & $2.37\times 10^{-3}$ \\
		MC & $5.74\times 10^{-3}$ & $4.55\times 10^{-3}$ \\
		PINN & $3.16\times 10^{-3}$ & $1.89\times 10^{-3}$ \\
		\bottomrule
	\end{tabular}\\
	\vspace{0.02\textwidth}
    \begin{minipage}[!h]{\textwidth}
        \centering
        $(b)$ Errors of $E$ at $t=0.5, 1.0$.
    \end{minipage}
    \vspace{-0.005\textwidth}
\end{minipage}
    %
    %
    %subTable3
  \begin{minipage}[!h]{\textwidth}
	\centering
% 	\makeatletter\def\@captype{}
	\begin{tabular}{ccc}
		\toprule
		$\ell^2(\text{energy})$ & up to $t=1.0$ & up to $t=5.0$ \\
		\midrule
		MM & $1.82\times 10^{-2}$  & $/$ \\
		MC & $/$ & $8.80\times 10^{-2}$ \\
		PINN & $3.60\times 10^{-2}$ & $/$ \\
		\bottomrule
	\end{tabular}
	\quad
	\begin{tabular}{ccc}
		\toprule
		$\text{RMSE}(\text{energy})$ & up to $t=1.0$ & up to $t=5.0$ \\
		\midrule
		MM & $3.97\times 10^{-3}$  & $/$ \\
		MC & $/$ & $9.33\times 10^{-3}$ \\
		PINN & $8.21\times 10^{-3}$ & $/$ \\
		\bottomrule
	\end{tabular}\\
	\vspace{0.02\textwidth}
    \begin{minipage}[!h]{\textwidth}
        \centering
        $(c)$ Errors of energy up to $t=1.0$ or $5.0$.
    \end{minipage}
    \vspace{-0.005\textwidth}
\end{minipage}
    \vspace{-0.04\textwidth}
\makeatletter\def\@captype{table}
	\caption{The linear Landau damping problem in kinetic regime $(\varepsilon=1.0)$. Relative $\ell^2$ error and RMSE of (a) density $\rho$, (b) electric field $E$ and (c) $\text{electric energy}$ for the micro-macro decomposition based APNN method (MM for short), the mass conservation based APNN method (MC for short) and vanilla PINN method (PINN for short).}
	\label{ld table kinetic}
\end{minipage}

%%%%%%%%%%%%%% ld t2 %%%%%%%%%%%%%%%%%%%%%

\begin{minipage}[htbp]{\textwidth}
	\centering
	%subTable1
\begin{minipage}[!h]{\textwidth}
	\centering
% 	\makeatletter\def\@captype{}
	\begin{tabular}{ccc}
		\toprule
		$\ell^2(\rho)$ & $t=0.5$ & $t=1.0$ \\
		\midrule
		MM & $9.46\times 10^{-4}$  & $1.70\times 10^{-3}$ \\
		MC & $3.50\times 10^{-3}$ & $2.18\times 10^{-3}$ \\
		PINN & $2.58\times 10^{-3}$ & $3.13\times 10^{-3}$ \\
		\bottomrule
	\end{tabular}
	\quad
	\begin{tabular}{ccc}
		\toprule
		$\text{RMSE}(\rho)$ & $t=0.5$ & $t=1.0$ \\
		\midrule
		MM & $9.46\times 10^{-4}$  & $1.70\times 10^{-3}$ \\
		MC & $3.50\times 10^{-3}$ & $2.18\times 10^{-3}$  \\
		PINN & $2.58\times 10^{-3}$ & $3.13\times 10^{-3}$ \\
		\bottomrule
	\end{tabular}\\
	\vspace{0.02\textwidth}
    \begin{minipage}[!h]{\textwidth}
        \centering
        $(a)$ Errors of $\rho$ at $t=0.5, 1.0$.
    \end{minipage}
    \vspace{-0.005\textwidth}
\end{minipage}
    %
    %
    %subtable2
\begin{minipage}[htbp]{\textwidth}
	\centering
		\begin{tabular}{ccc}
		\toprule
		$\ell^2(E)$ & $t=0.5$ & $t=1.0$ \\
		\midrule
		MM & $2.00\times 10^{-2}$  & $6.35\times 10^{-2}$ \\
		MC & $1.25\times 10^{-1}$ & $1.34\times 10^{-1}$ \\
		PINN & $9.09\times 10^{-2}$ & $1.72\times 10^{-1}$ \\
		\bottomrule
	\end{tabular}
	\quad
	\begin{tabular}{ccc}
		\toprule
		$\text{RMSE}(E)$ & $t=0.5$ & $t=1.0$ \\
		\midrule
		MM & $1.13\times 10^{-3}$  & $2.07\times 10^{-3}$ \\
		MC & $7.03\times 10^{-3}$ & $4.37\times 10^{-3}$ \\
		PINN & $5.13\times 10^{-3}$ & $5.59\times 10^{-3}$ \\
		\bottomrule
	\end{tabular}\\
	\vspace{0.02\textwidth}
    \begin{minipage}[!h]{\textwidth}
        \centering
        $(b)$ Errors of $E$ at $t=0.5, 1.0$.
    \end{minipage}
    \vspace{-0.005\textwidth}
\end{minipage}
    %
    %
    %subTable3
  \begin{minipage}[!h]{\textwidth}
	\centering
% 	\makeatletter\def\@captype{}
	\begin{tabular}{ccc}
		\toprule
		$\ell^2(\text{energy})$ & up to $t=1.0$ & up to $t=5.0$ \\
		\midrule
		MM & $1.23\times 10^{-2}$  & $/$ \\
		MC & $/$ & $1.05\times 10^{-1}$ \\
		PINN & $7.22\times 10^{-2}$ & $/$ \\
		\bottomrule
	\end{tabular}
	\quad
	\begin{tabular}{ccc}
		\toprule
		$\text{RMSE}(\text{energy})$ & up to $t=1.0$ & up to $t=5.0$ \\
		\midrule
		MM & $2.53\times 10^{-3}$  & $/$ \\
		MC & $/$ & $1.01\times 10^{-2}$ \\
		PINN & $1.49\times 10^{-2}$ & $/$ \\
		\bottomrule
	\end{tabular}\\
	\vspace{0.02\textwidth}
    \begin{minipage}[!h]{\textwidth}
        \centering
        $(c)$ Errors of energy up to $t=1.0$ or $5.0$.
    \end{minipage}
    \vspace{-0.005\textwidth}
\end{minipage}
    \vspace{-0.04\textwidth}
\makeatletter\def\@captype{table}
	\caption{Errors of linear Landau damping with intermediate regime $(\varepsilon=0.5)$. Relative $\ell^2$ error and RMSE of (a) density $\rho$, (b) electric field $E$ and (c) $\text{electric energy}$ for the micro-macro decomposition based APNN method (MM for short), the mass conservation based APNN method (MC for short) and vanilla PINN method (PINN for short).}
	\label{ld table inter}
\end{minipage}

%%%%%%%%%%%%%% ld t3 %%%%%%%%%%%%%%%%%%%%%

\begin{minipage}[htbp]{\textwidth}
	\centering
	%subTable1
\begin{minipage}[!h]{\textwidth}
	\centering
% 	\makeatletter\def\@captype{}
	\begin{tabular}{ccc}
		\toprule
		$\ell^2(\rho)$ & $t=0.5$ & $t=1.0$ \\
		\midrule
		MM & $3.80\times 10^{-4}$  & $5.92\times 10^{-4}$ \\
		MC & $1.12\times 10^{-3}$ & $8.00\times 10^{-4}$ \\
		PINN & $1.03\times 10^{-2}$ & $1.84\times 10^{-2}$ \\
		\bottomrule
	\end{tabular}
	\quad
	\begin{tabular}{ccc}
		\toprule
		$\text{RMSE}(\rho)$ & $t=0.5$ & $t=1.0$ \\
		\midrule
		MM & $3.80\times 10^{-4}$  & $5.92\times 10^{-4}$ \\
		MC & $1.12\times 10^{-3}$ & $8.00\times 10^{-4}$  \\
		PINN & $1.03\times 10^{-2}$ & $1.84\times 10^{-2}$ \\
		\bottomrule
	\end{tabular}\\
	\vspace{0.02\textwidth}
    \begin{minipage}[!h]{\textwidth}
        \centering
        $(a)$ Errors of $\rho$ at $t=0.5, 1.0$.
    \end{minipage}
    \vspace{-0.005\textwidth}
\end{minipage}
    %
    %
    %subtable2
\begin{minipage}[htbp]{\textwidth}
	\centering
		\begin{tabular}{ccc}
		\toprule
		$\ell^2(E)$ & $t=0.5$ & $t=1.0$ \\
		\midrule
		MM & $8.73\times 10^{-3}$  & $3.74\times 10^{-2}$ \\
		MC & $5.00\times 10^{-2}$ & $6.10\times 10^{-2}$ \\
		PINN & $4.78\times 10^{-1}$ & $1.42$ \\
		\bottomrule
	\end{tabular}
	\quad
	\begin{tabular}{ccc}
		\toprule
		$\text{RMSE}(E)$ & $t=0.5$ & $t=1.0$ \\
		\midrule
		MM & $3.73\times 10^{-4}$  & $9.61\times 10^{-4}$ \\
		MC & $2.13\times 10^{-3}$ & $1.57\times 10^{-3}$ \\
		PINN & $2.04\times 10^{-2}$ & $3.65\times 10^{-2}$ \\
		\bottomrule
	\end{tabular}\\
	\vspace{0.02\textwidth}
    \begin{minipage}[!h]{\textwidth}
        \centering
        $(b)$ Errors of $E$ at $t=0.5, 1.0$.
    \end{minipage}
    \vspace{-0.005\textwidth}
\end{minipage}
    %
    %
    %subTable3
  \begin{minipage}[!h]{\textwidth}
	\centering
% 	\makeatletter\def\@captype{}
	\begin{tabular}{ccc}
		\toprule
		$\ell^2(\text{energy})$ & up to $t=1.0$ & up to $t=5.0$ \\
		\midrule
		MM & $7.03\times 10^{-2}$  & $/$ \\
		MC & $/$ & $3.67\times 10^{-2}$ \\
		PINN & $4.56\times 10^{-1}$ & $/$ \\
		\bottomrule
	\end{tabular}
	\quad
	\begin{tabular}{ccc}
		\toprule
		$\text{RMSE}(\text{energy})$ & up to $t=1.0$ & up to $t=5.0$ \\
		\midrule
		MM & $5.77\times 10^{-3}$  & $/$ \\
		MC & $/$ & $3.01\times 10^{-3}$ \\
		PINN & $8.10\times 10^{-2}$ & $/$ \\
		\bottomrule
	\end{tabular}\\
	\vspace{0.02\textwidth}
    \begin{minipage}[!h]{\textwidth}
        \centering
        $(c)$ Errors of energy up to $t=1.0$ or $5.0$.
    \end{minipage}
    \vspace{-0.005\textwidth}
\end{minipage}
    \vspace{-0.04\textwidth}
\makeatletter\def\@captype{table}
	\caption{Errors of linear Landau damping with high-field regime $(\varepsilon=0.01)$. Relative $\ell^2$ error and RMSE of (a) density $\rho$, (b) electric field $E$ and (c) $\text{electric energy}$ for the micro-macro decomposition based APNN method (MM for short), the mass conservation based APNN method (MC for short) and vanilla PINN method (PINN for short).}
	\label{ld table high}
\end{minipage}

%%%%%%%%%%%%%%%%%%%% DP %%%%%%%%%%%%%%%%%%%%%
\subsection{Problem II: The bump-on-tail case}

In this part, we test whether proposed APNN methods are still effective when initial data is non-equilibrium. Consider the “double peak” function initial data:
\begin{equation*}
    f_0(x, v)=\frac{1}{\sqrt{2 \pi}}\left(\frac{9}{10} \exp \left(-\frac{v^2}{2}\right)+\frac{2}{10} \exp \left(-4(v-4.5)^2\right)\right)(1+\alpha \cos (k x)),
\end{equation*}
where $\alpha=0.05, k=0.5$ and $(x, v) \in[0,2 \pi / k] \times\left[-v_{\max }, v_{\max }\right]$ with $v_{\max }=8$. The  ICs for $\rho, \phi$ are the same as the Landau damping case. The periodic BCs are applied.

\begin{figure}[htbp]
    \centering 
    \includegraphics[width=12cm]{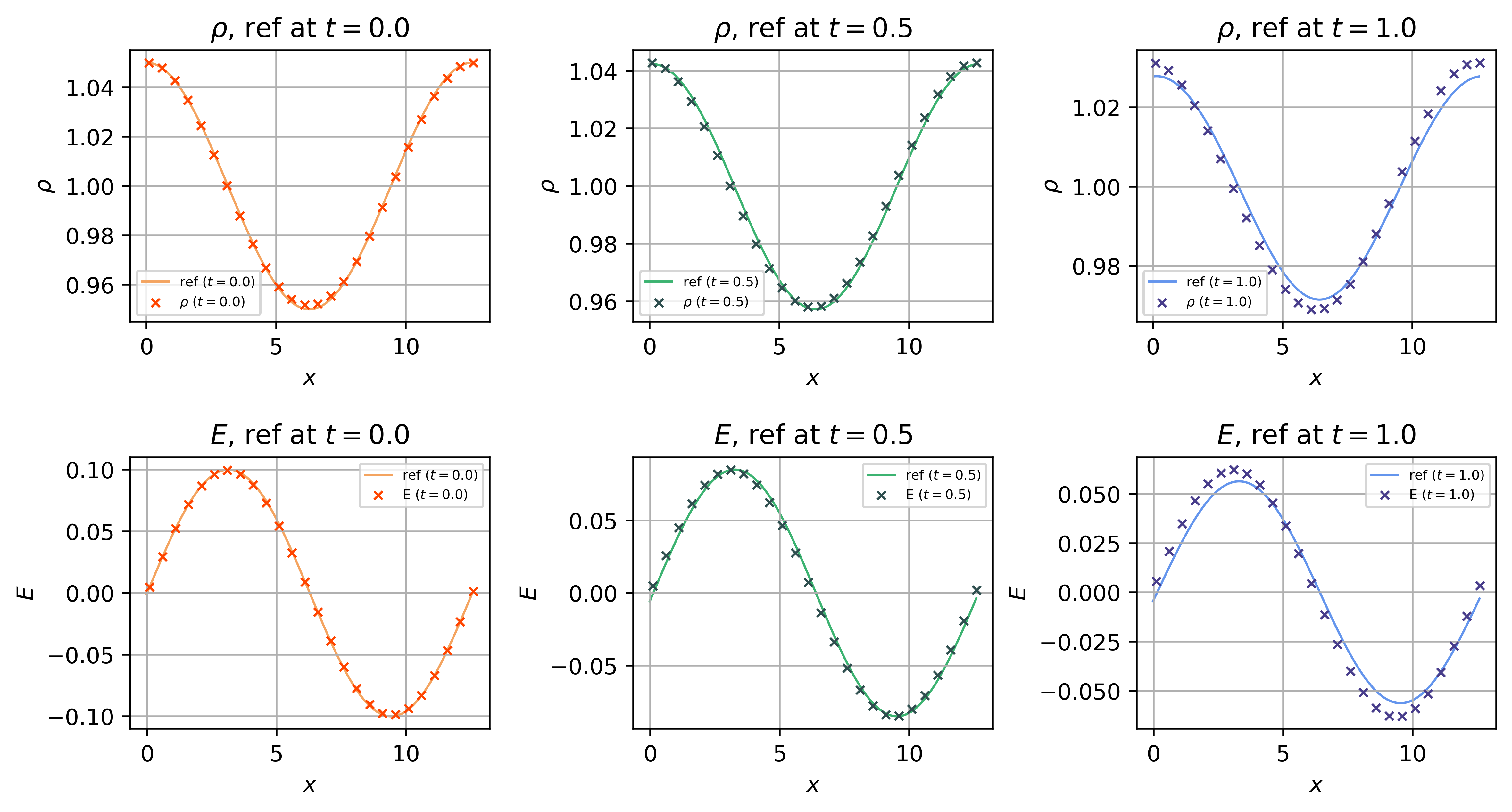}
    \caption{The bump-on-tail case solved by the micro-macro decomposition based APNN method with kinetic regime ($\varepsilon=1.0)$. Density $\rho$ as a function of space $x$ for $t = 0.0, 0.5, 1.0$ (top). Electric field $E$ as a function of space $x$ for $t = 0.0, 0.5, 1.0$ (bottom). Neural networks are $[3, 128, 128, 128, 128, 128, 1]$ for $\rho, \phi$ and $[4, 256, 256, 256, 256, 256, 1]$ for $g$. Batch size is $512$ in domain and $256$ for initial condition. Penalty $\lambda_1 = 50$ and $\lambda_3 = 1$. Errors: Relative $\ell^2$ error of $\rho$ is $1.93\times 10^{-3}$ at $t=0.5$ and $3.01\times 10^{-3}$ at $t=1.0$. Relative $\ell^2$ error of $E$ is $6.28\times 10^{-2}$ at $t=0.5$ and $1.58\times 10^{-1}$ at $t=1.0$.}
    \label{dp mm kinetic}
\end{figure}

\begin{figure}[htbp]
    \centering 
    \includegraphics[width=12cm]{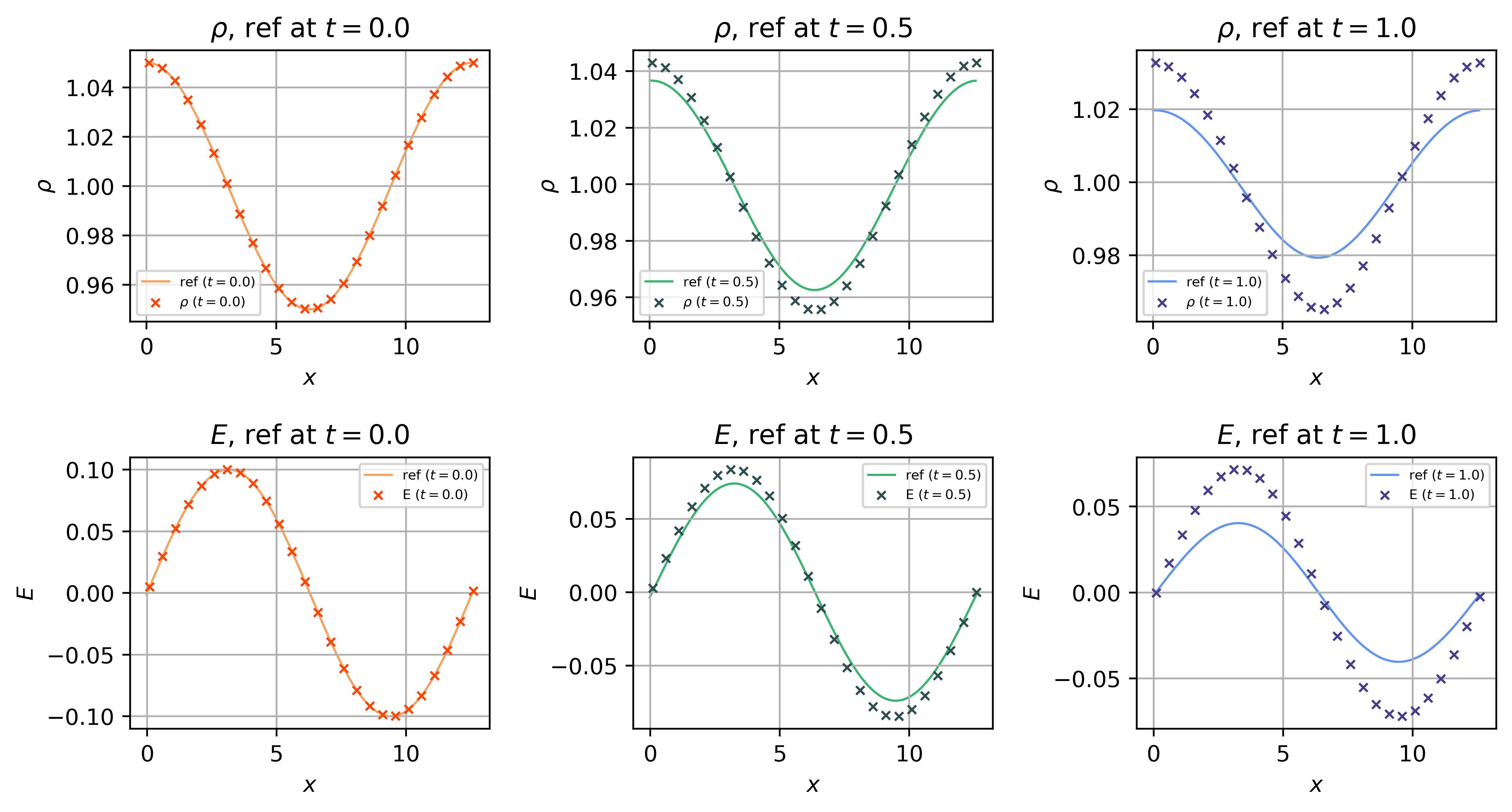}
    \caption{The bump-on-tail case solved by the micro-macro decomposition based APNN method with intermediate regime ($\varepsilon=0.3)$. Density $\rho$ as a function of space $x$ for $t = 0.0, 0.5, 1.0$ (top). Electric field $E$ as a function of space $x$ for $t = 0.0, 0.5, 1.0$ (bottom). Neural networks and batch size are the same as kinetic regime case. Penalty $\lambda_1 = 1$ and $\lambda_3 = 1000$. Errors: Relative $\ell^2$ error of $\rho$ is $4.76\times 10^{-3}$ at $t=0.5$ and $9.58\times 10^{-3}$ at $t=1.0$. Relative $\ell^2$ error of $E$ is $1.44\times 10^{-1}$ at $t=0.5$ and $7.86\times 10^{-1}$ at $t=1.0$.}
    \label{dp mm inter}
\end{figure}

\begin{figure}[htbp]
    \centering 
    \includegraphics[width=12cm]{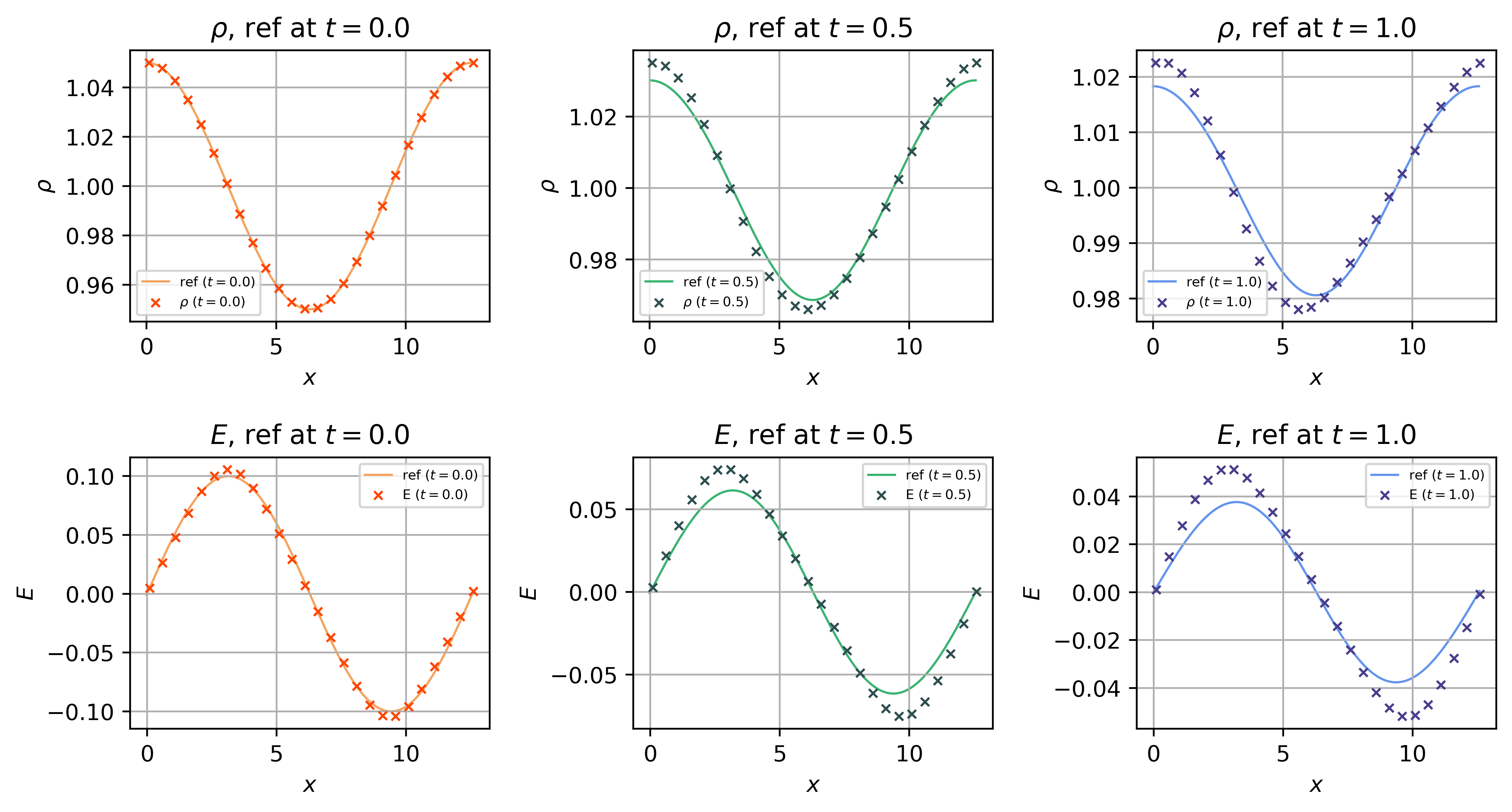}
    \caption{The bump-on-tail case solved by the micro-macro decomposition based APNN method with high-field regime ($\varepsilon=0.001)$. Density $\rho$ as a function of space $x$ for $t = 0.0, 0.5, 1.0$ (top). Electric field $E$ as a function of space $x$ for $t = 0.0, 0.5, 1.0$ (bottom). Neural networks, batch size and penalty are the same as intermediate regime case. Errors: Relative $\ell^2$ error of $\rho$ is $2.94\times 10^{-3}$ at $t=0.5$ and $2.98\times 10^{-3}$ at $t=1.0$. Relative $\ell^2$ error of $E$ is $2.10\times 10^{-1}$ at $t=0.5$ and $3.63\times 10^{-1}$ at $t=1.0$.}
    \label{dp mm high}
\end{figure}

\begin{figure}[htbp]
	\centering  %图片全局居中
	\subfigbottomskip=4pt %两行子图之间的行间距
	\subfigcapskip=-8pt %设置子图与子标题之间的距离
	\subfigure[Kinetic regime $(\varepsilon=1)$.]{
		\includegraphics[width=0.6\linewidth]{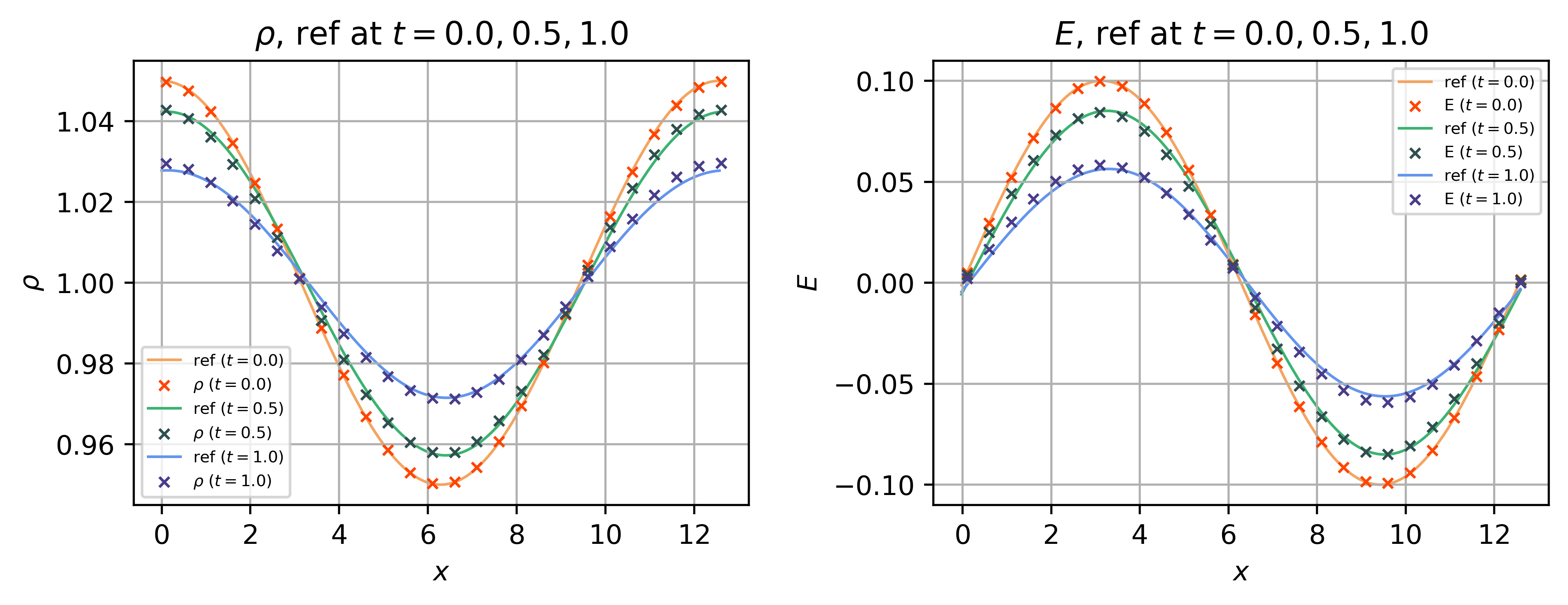}}
	\subfigure[Intermediate regime $(\varepsilon=0.3)$.]{
		\includegraphics[width=0.6\linewidth]{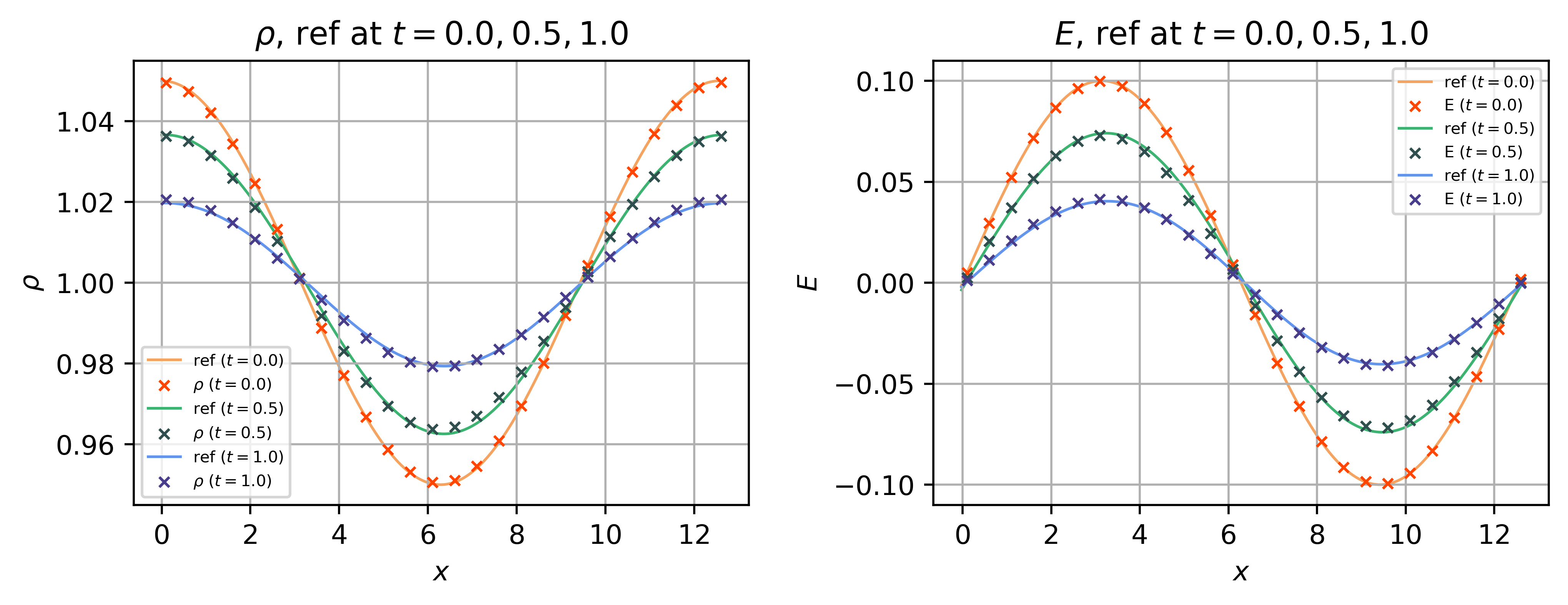}}
	\subfigure[High-field regime $(\varepsilon=0.001)$.]{
		\includegraphics[width=0.6\linewidth]{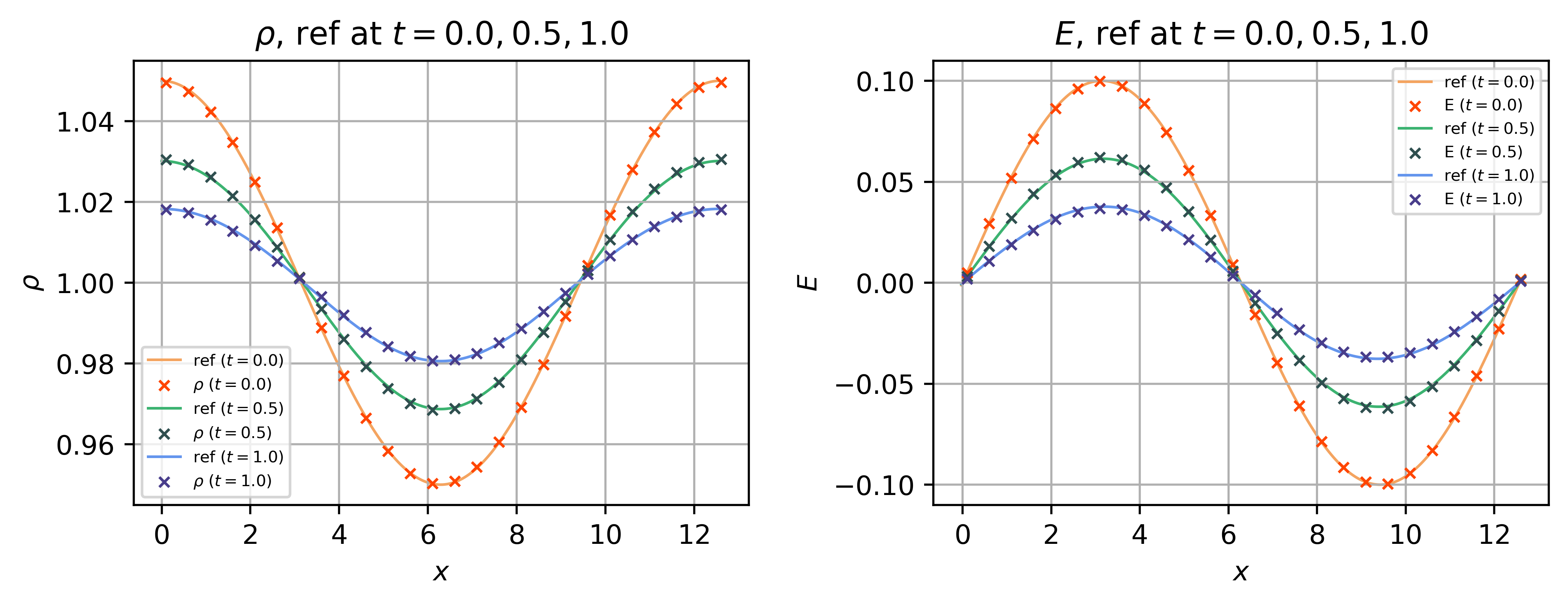}}
	\caption{The bump-on-tail case solved by the mass conservation based APNN method. Density $\rho$ (left column) and electric field $E$ (right column) as functions of space $x$ at $t = 0.0, 0.5, 1.0$. Neural networks are $[3, 128, 128, 128, 128, 128, 1]$ for $\rho, \phi$ and $[4, 256, 256, 256, 256, 256, 1]$ for $f$. Batch size is $512$ in domain, $256$ for initial condition and $256$ for conservation condition $\rho=\langle f\rangle$. Penalty $\kappa_1 = 50$, $\kappa_3 = \kappa_4 = 1$ for each regimes. Errors: $(a)$ Relative $\ell^2$ error of $\rho$ is $1.51\times 10^{-3}$ at $t=0.5$ and $1.32\times 10^{-3}$ at $t=1.0$. Relative $\ell^2$ error of $E$ is $4.90\times 10^{-2}$ at $t=0.5$ and $6.53\times 10^{-2}$ at $t=1.0$. $(b)$ Relative $\ell^2$ error of $\rho$ is $1.02\times 10^{-3}$ at $t=0.5$ and $5.83\times 10^{-4}$ at $t=1.0$. Relative $\ell^2$ error of $E$ is $3.67\times 10^{-2}$ at $t=0.5$ and $3.22\times 10^{-2}$ at $t=1.0$. $(c)$ Relative $\ell^2$ error of $\rho$ is $1.93\times 10^{-3}$ at $t=0.5$ and $3.01\times 10^{-3}$ at $t=1.0$. Relative $\ell^2$ error of $E$ is $6.28\times 10^{-2}$ at $t=0.5$ and $1.58\times 10^{-1}$ at $t=1.0$.}
    \label{dp mc}
\end{figure}

Here, we demonstrate the performance of the two proposed models from kinetic regimes to high-field regimes. The micro-macro decomposition-based APNN method exhibits inaccuracy across all regimes, as depicted in Figure \ref{dp mm kinetic}, Figure \ref{dp mm inter} and Figure \ref{dp mm high}. This inaccuracy may be attributable to the severe variations in the magnitude of $g$, which occurs when the dynamics start from non-equilibrium initial data. As time progresses, these variation lead to a drastic change in order of $g$, adversely impacting the performance of the neural network function $g$. Consequently, this results in the inaccuracy of the learning outcomes \cite{lou2021physics}. Same as the micro-macro decomposition based APNN method for the Landau damping case, the inaccuracy of the learning outcomes attribute to the severe variations in the magnitude of $g$. In contrast, the mass conservation based APNN method maintains the accurate outputs of neural networks for $\rho$ and $f$ at a consistent order of $O(1)$, thereby ensuring good learning performance shown in Figure \ref{dp mc}. As illustrated by this case, the APNN method exhibits adaptability to an non-equilibrium initial data, which underlines its flexibility and robustness.

%%%%%%%%%%%%%%%%%%%% Riemann %%%%%%%%%%%%%%%%%%%%%
\subsection{Problem III: A Riemann problem}

In this part, a one-dimensional Riemann problem is used to examine the performance of two proposed APNN methods in discontinuity cases. The piecewise constant initial condition is defined as:
\begin{equation}
    \begin{cases}
        \left(\rho_l, h_l\right)=(1 / 8,1 / 2), & 0 \leq x<1 / 4, \\ 
        \left(\rho_m, h_m\right)=(1 / 2,1 / 8), & 1 / 4 \leq x<3/4,\\ 
        \left(\rho_r, h_r\right)=(1 / 8,1 / 2), & 3 / 4 \leq x \leq 1.
    \end{cases}
\end{equation}
The exact initial condition of $\phi$ is solved by the Poisson equation \ref{poisson} as:
\begin{equation}
    \begin{cases}
        \phi_l=3/16x^2-3/256, & 0 \leq x<1 / 4,\\ 
        \phi_m=-3/16(x-1/2)^2+3/256, & 1 / 4 \leq x<3 / 4,\\ 
        \phi_r=3/16(x-1)^2-3/256, & 3 / 4 \leq x \leq 1.
    \end{cases}
\end{equation}
The corresponding first-order derivative $\partial_x\phi$ is directly derived. Set $f$ initially as:
\begin{equation}
    f_0(x,v)=\frac{\rho_0}{\sqrt{2 \pi}}\exp\left(-\frac{(v+\partial_x\phi_0(x))^2}{2}\right). 
\end{equation}
Again periodic BC in $x$ direction is applied. 

\begin{figure}[htbp]
    \centering 
    \includegraphics[width=0.8\linewidth]{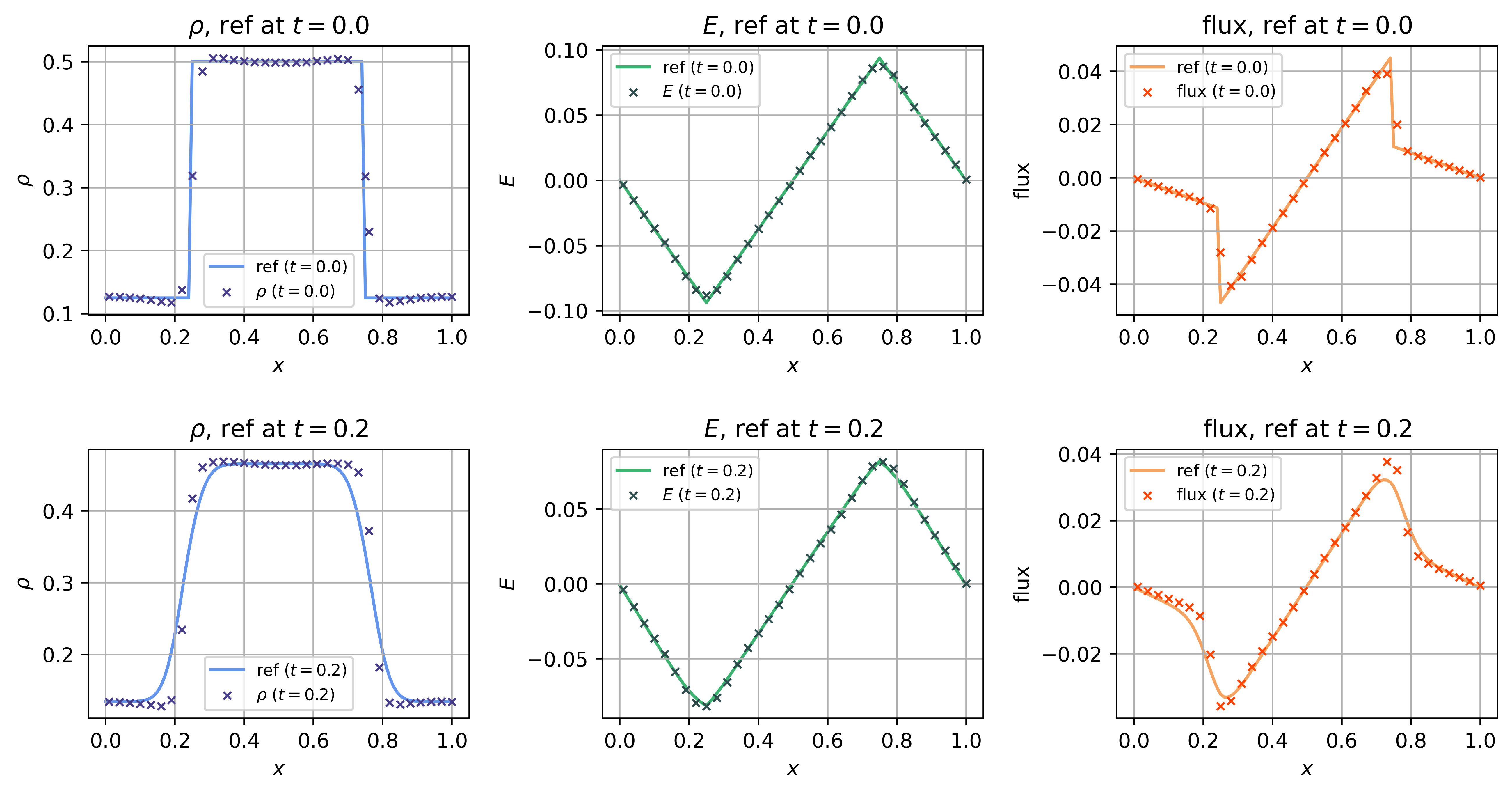}
    \caption{The Riemann problem by the micro-macro decomposition based APNN method in the high-field regime ($\varepsilon=0.001$). Density $\rho$ as a function of space $x$ at $t = 0.0, 0.2$ (left column). Electric field $E$ as a function of space $x$ at $t = 0.0, 0.2$ (middle column). Flux as a function of space $x$ at $t = 0.0, 0.2$ (right column). Neural networks are $[3, 128, 128, 128, 128, 128, 1]$ for $\rho, \phi$ and $[4, 256, 256, 256, 256, 256, 1]$ for $g$. Batch size is $512$ in domain and $512$ for initial condition. Penalty $\lambda_1^{\text{macro}} = \lambda_1^{\text{micro}} = \lambda_3= 5$ and the rest of penalties equal to $1$.}
    \label{riemann mm}
\end{figure}

\begin{figure}[htbp]
    \centering 
    \includegraphics[width=0.8\linewidth]{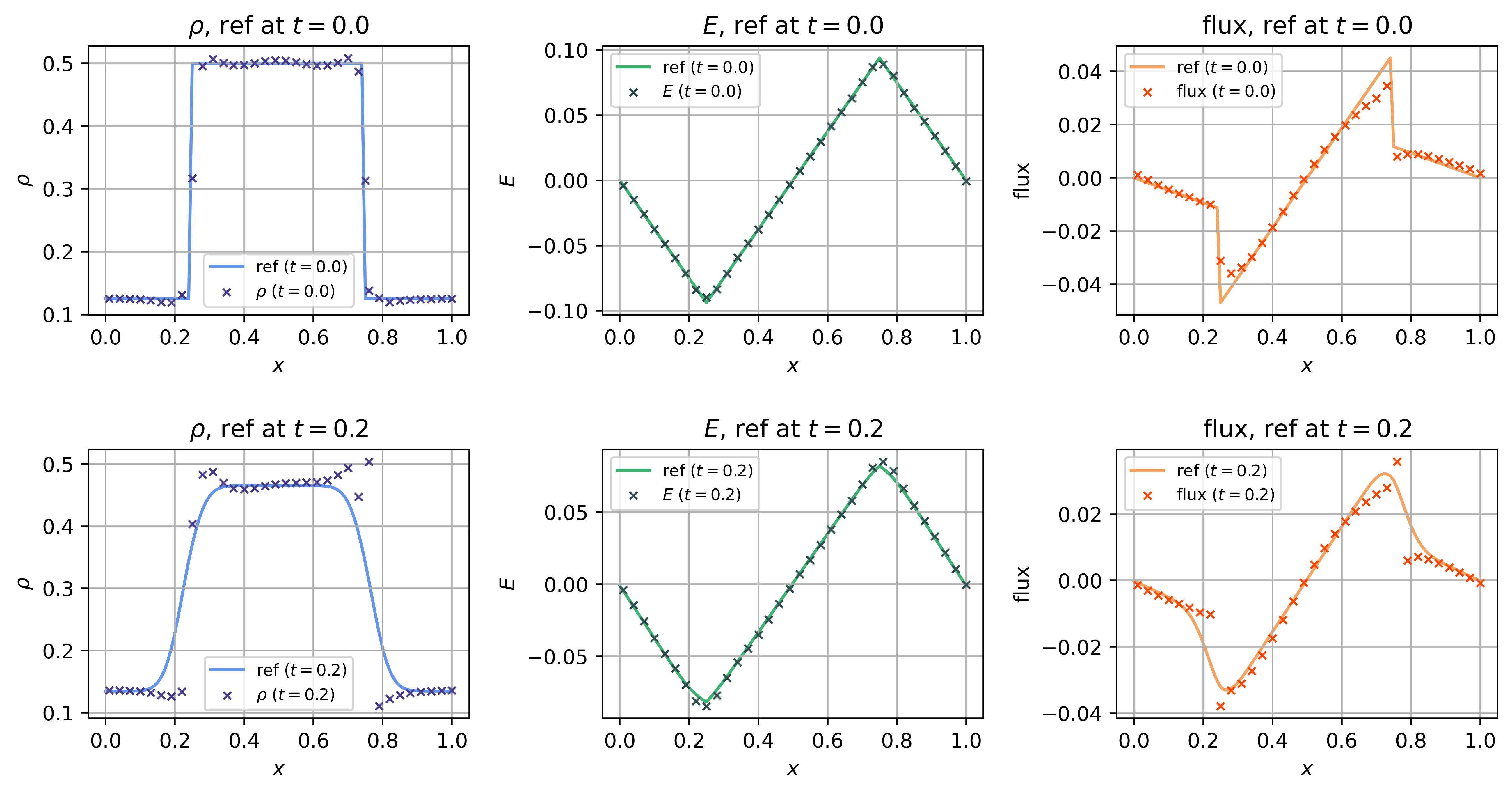}
    \caption{The Riemann problem by the mass conservation based APNN method in the high-field regime ($\varepsilon=0.001$). Density $\rho$ as a function of space $x$ at $t = 0.0, 0.2$ (left column). Electric field $E$ as a function of space $x$ at $t = 0.0, 0.2$ (middle column). Flux as a function of space $x$ at $t = 0.0, 0.2$ (right column). Neural networks are $[3, 128, 128, 128, 128, 128, 1]$ for $\rho, \phi$ and $[4, 256, 256, 256, 256, 256, 1]$ for $f$. Batch size is $512$ in domain, $512$ for initial condition and $256$ for conservation condition $\rho=\langle f\rangle$. Penalty $\kappa_3^\rho = 3$, $\kappa_3^f = 1000$ and the rest of penalties equal to $1$.}
    \label{riemann mc}
\end{figure}

We check the ability of two proposed APNN methods in tackling the discontinuities. Macroscopic variables $\rho, \phi$ and flux $j(t, x)=\int_{\mathcal{R}} v f(t, x, v) \mathrm{d} v$ are plotted in the high-field case $(\varepsilon=0.001)$. Figure \ref{riemann mm} and Figure \ref{riemann mc} show that the proposed APNN methods can capture the sharp transition at singularities of the solution. High-frequency information such as sharp transitions are challenging to be learned by neural networks \cite{xu2019frequency}. This difficulty can be covered by both APNN methods in this case.

We further note that the micro-macro decomposition based APNN method provides more accurate profiles in a shorter time frame, compared to the mass conservation based APNN method. This observation is supported by Table \ref{riemann table} and the bottom rows of Figure \ref{riemann mm} and Figure \ref{riemann mc}. One possible explanation is due to the explicit representation of the non-equilibrium component $g$ in the micro-macro decomposition based APNN method.This explicit representation enhances the DNNs' ability to accurately capture the non-equilibrium part of the distribution. The non-equilibrium part has a strong influence on the convergence of the governing VPFP system \cite{lou2021physics}. As a result, DNN could provide more detailed and insightful information. However, it's important to highlight that this improvement in accuracy is not universal. Specifically, it is only valid when the calculated time is relatively short and $\rho$ and $g$ are comparable in magnitude.

\begin{minipage}[htbp]{\textwidth}
	\centering
	%subTable1
\begin{minipage}[!h]{\textwidth}
	\centering
% 	\makeatletter\def\@captype{}
	\begin{tabular}{ccc}
		\toprule
		$\ell^2(\rho)$ & $t=0.0$ & $t=0.2$ \\
		\midrule
		MM & $9.55\times 10^{-2}$  & $7.60\times 10^{-2}$     \\
		MC & $7.36\times 10^{-2}$ & $1.34\times 10^{-2}$      \\
		\bottomrule
	\end{tabular}
	\quad
	\begin{tabular}{ccc}
		\toprule
		$\text{RMSE}(\rho)$ & $t=0.0$ & $t=0.2$ \\
		\midrule
		MM & $3.46\times 10^{-2}$  & $2.61\times 10^{-2}$     \\
		MC & $2.67\times 10^{-2}$ & $4.59\times 10^{-2}$     \\
		\bottomrule
	\end{tabular}\\
	\vspace{0.02\textwidth}
    \begin{minipage}[!h]{\textwidth}
        \centering
        $(a)$ Errors of $\rho$ at $t=0.0, 0.2$.
    \end{minipage}
    \vspace{-0.005\textwidth}
\end{minipage}
    %
    %
    %subtable2
\begin{minipage}[htbp]{\textwidth}
	\centering
		\begin{tabular}{ccc}
		\toprule
		$\ell^2(E)$ & $t=0.0$ & $t=0.2$ \\
		\midrule
		MM & $1.61\times 10^{-2}$  & $3.18\times 10^{-2}$ \\
		MC & $1.54\times 10^{-2}$ & $4.06\times 10^{-2}$ \\
		\bottomrule
	\end{tabular}
	\quad
	\begin{tabular}{ccc}
		\toprule
		$\text{RMSE}(E)$ & $t=0.0$ & $t=0.2$ \\
		\midrule
		MM & $1.40\times 10^{-3}$  & $1.58\times 10^{-3}$ \\
		MC & $8.32\times 10^{-4}$ & $2.01\times 10^{-3}$ \\
		\bottomrule
	\end{tabular}\\
	\vspace{0.02\textwidth}
    \begin{minipage}[!h]{\textwidth}
        \centering
        $(b)$ Errors of $E$ at $t=0.0, 0.2$.
    \end{minipage}
    \vspace{-0.005\textwidth}
\end{minipage}
    %
    %
    %subTable3
  \begin{minipage}[!h]{\textwidth}
	\centering
% 	\makeatletter\def\@captype{}
	\begin{tabular}{ccc}
		\toprule
		$\ell^2(\text{flux})$ & $t=0.0$ & $t=0.2$ \\
		\midrule
		MM & $1.64\times 10^{-1}$  & $1.46\times 10^{-1}$ \\
		MC & $1.76\times 10^{-1}$ & $2.34\times 10^{-1}$  \\
		\bottomrule
	\end{tabular}
	\quad
	\begin{tabular}{ccc}
		\toprule
		$\text{RMSE}(\text{flux})$ & $t=0.0$ & $t=0.2$ \\
		\midrule
		MM & $3.23\times 10^{-3}$  & $2.60\times 10^{-3}$  \\
		MC & $3.46\times 10^{-3}$ & $4.17\times 10^{-3}$ \\
		\bottomrule
	\end{tabular} \\
	\vspace{0.02\textwidth}
    \begin{minipage}[!h]{\textwidth}
        \centering
        $(c)$ Errors of flux at $t=0.0, 0.2$.
    \end{minipage}
    \vspace{-0.005\textwidth}
\end{minipage}
    \vspace{-0.04\textwidth}
\makeatletter\def\@captype{table}
	\caption{Errors of Riemann problem in the high-field regime $(\varepsilon=0.001)$. Relative $\ell^2$ error and RMSE of $\rho, E, \text{flux}$ at $t=0.0, 0.2$  for the micro-macro decomposition based APNN method (MM for short) and the mass conservation based APNN method (MC for short).}
	\label{riemann table}
\end{minipage}

%%%%%%%%%%%%%%%%%%%% Mixing %%%%%%%%%%%%%%%%%%%%%
\subsection{Problem IV: Mixing regimes}

The proposed APNN methods are now tested in the case with mixing regimes. In these regimes, the scale parameter $\varepsilon$ fluctuates across several orders of magnitude in space. The parameter $\varepsilon$ is constructed as:
\begin{equation}
    \varepsilon(x)= 
    \begin{cases}
        \varepsilon_0+\frac{1}{2}\left(\tanh (5-10 x)+\tanh (5+10 x)\right), & -1 \leq x \leq 0.3,\\ 
        \varepsilon_0, & 0.3 <x \leq 1,
    \end{cases}
\end{equation}
where $\varepsilon_0$ is set as $0.001$. Hence both the kinetic and high-field regimes are under consideration. The initial conditions for $\rho, f$ are given by:
\begin{equation}
    \rho_0(x)=\frac{\sqrt{2 \pi}}{6}(2+\sin (\pi x)), \quad f_0(x,v)=\frac{\rho_0(x)}{\sqrt{2 \pi}}\exp(-\frac{\left(v+\partial_x\phi_0(x)\right)^2}{2}).
\end{equation}
The initial electric field $\phi_0(x)$ is analytically solved by the Poisson equation \ref{poisson} with $h(x)=\sqrt{2\pi}/3$. Periodic BC in $x$ direction is imposed. 

\begin{figure}[htbp]
    \centering 
    \includegraphics[width=5.0cm]{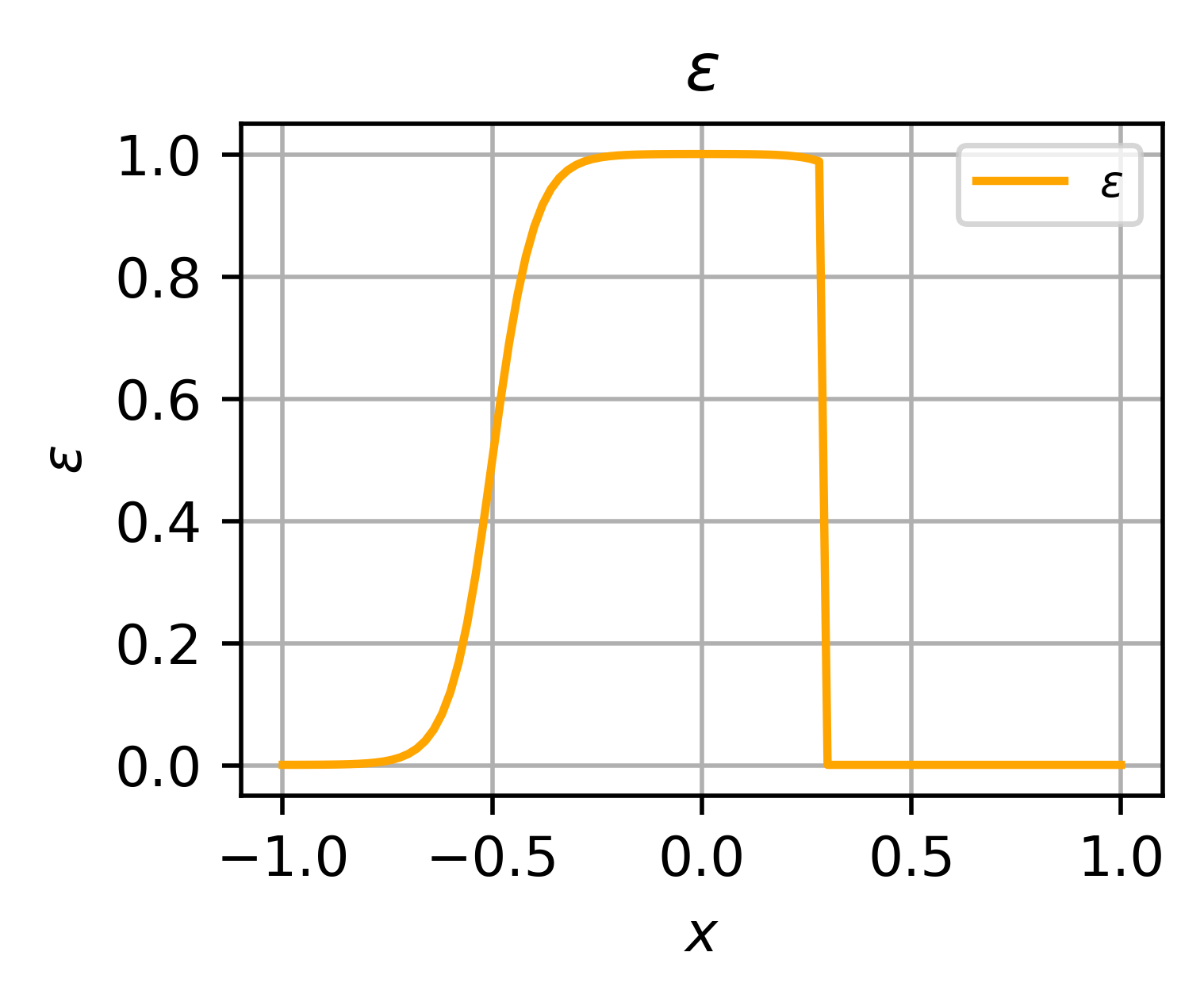}
    \caption{The scale parameter $\varepsilon$ as a function of space $x$ with several orders of magnitude.}
    \label{mixing eps}
\end{figure}

\begin{figure}[htbp]
	\centering  %图片全局居中
	\subfigbottomskip=6pt %两行子图之间的行间距
	\subfigcapskip=-8pt %设置子图与子标题之间的距离
	\subfigure[The micro-macro decomposition based APNN method.]{
		\includegraphics[width=0.8\linewidth]{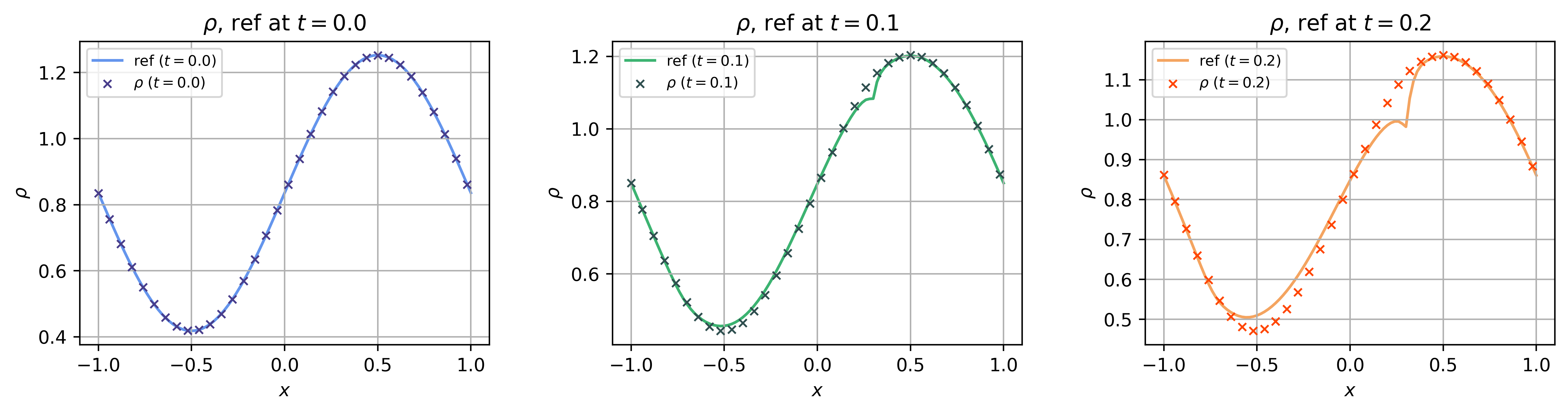}}
	\subfigure[The mass conservation based APNN method.]{
		\includegraphics[width=0.8\linewidth]{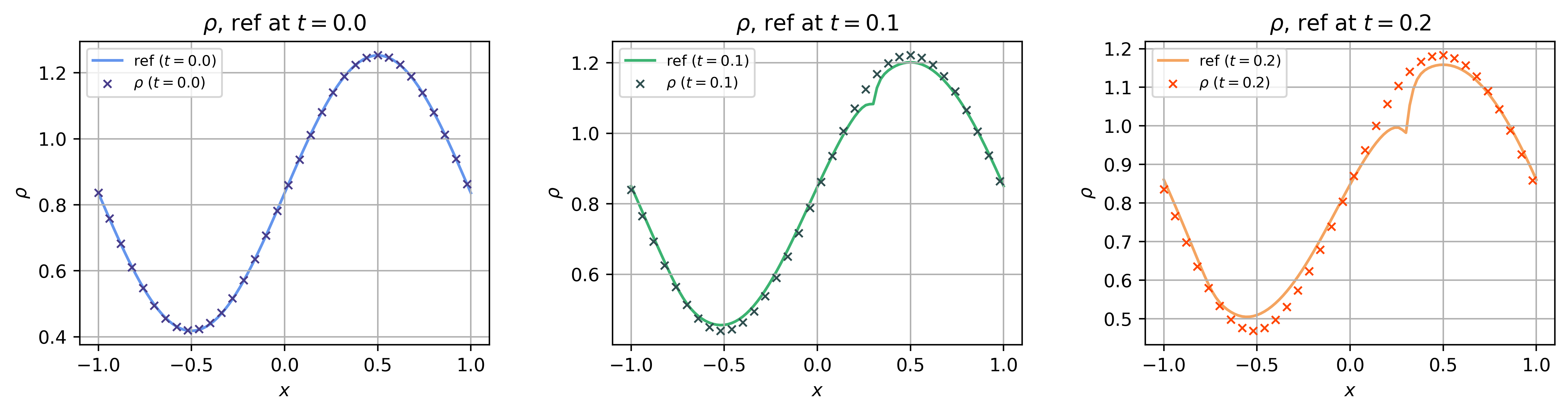}}
	\caption{The mixing regimes problem solved by the micro-macro decomposition based and mass conservation based APNN methods. Density $\rho$ as a function of space $x$ at $t = 0.0, 0.1, 0.2$. $(a)$ Neural networks are $[3, 128, 128, 128, 128, 128, 1]$ for $\rho, \phi$ and $[4, 256, 256, 256, 256, 256, 1]$ for $g$. Batch size is $512$ in domain and $256$ for initial condition.  Penalty $\lambda_1 = 0.1$ and $\lambda_3 = 1$. $(b)$ Neural networks are $[3, 128, 128, 128, 128, 128, 1]$ for $\rho, \phi$ and $[4, 256, 256, 256, 256, 256, 1]$ for $f$. Batch size is $512$ in domain, $256$ for initial condition and $256$ for conservation condition $\rho=\langle f\rangle$. Penalty $\kappa_1 = \kappa_3 = \kappa_4 = 1$.}
    \label{mixing}
\end{figure}

We report the capacity of the proposed APNN methods to resolve the case with mixing regimes. Figure \ref{mixing eps} plots the scale parameter $\varepsilon$ with a discontinuity at $x=0.3$. The small $\varepsilon$ induces transition of the solution is hard to be learned by DNNs-based methods as explained in the Riemann problem. Figure \ref{mixing} and Table \ref{mixing table} show that both proposed APNN methods capture quite well the profile of density $\rho$, except at $x=0.3$ where the solution has a weak discontinuity. This inferior performance  is ascribed to the inherent difficulty in learning high-frequency information by neural networks.

\begin{minipage}[htbp]{\textwidth}
	\makeatletter\def\@captype{table}
	\centering
	\begin{tabular}{cccc}
		\toprule
		$\ell^2(\rho)$ & $t=0.0$ & $t=0.1$ & $t=0.2$ \\
		\midrule
		MM & $9.23\times 10^{-4}$ &  $1.36\times 10^{-2}$ & $3.67\times 10^{-2}$  \\
		MC & $1.85\times 10^{-3}$ &  $2.06\times 10^{-2}$ & $4.40\times 10^{-2}$  \\
		\bottomrule
	\end{tabular}
	\caption{Error of the mixing regimes problem. Relative $\ell^2$ error of $\rho$ at $t=0.0, 0.1, 0.2$ for the micro-macro decomposition based method (MM for short) and the mass conservation based method (MC for short).}
	\label{mixing table}
\end{minipage}

%%%%%%%%%%%%%%%%%%%% Gravity %%%%%%%%%%%%%%%%%%%%%
\subsection{Problem V: The gravitational case}
 This section illustrates the efficacy of the proposed APNN methods for the VPFP system with gravitational force. Compared to the previous cases under repulsive electrostatic interaction, the attractive gravitational force induces different behavior. The existence of a unique weak solution locally in time has been proved in \cite{nieto2001high} for the limit equation with gravitational force.
 
\begin{figure}[htbp]
	\centering  %图片全局居中
	\subfigbottomskip=6pt %两行子图之间的行间距
	\subfigcapskip=-8pt %设置子图与子标题之间的距离
	\subfigure[The micro-macro decomposition based APNN method.]{
		\includegraphics[width=0.6\linewidth]{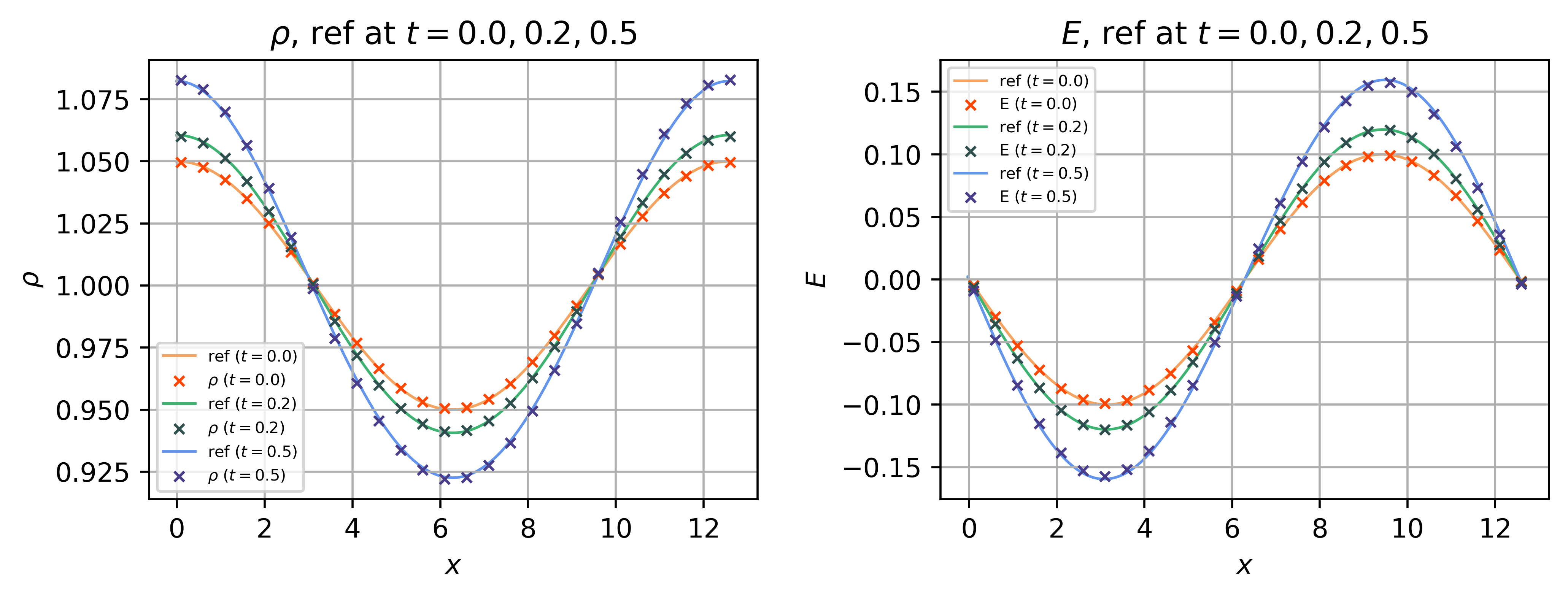}}
	\subfigure[The mass conservation based APNN method.]{
		\includegraphics[width=0.6\linewidth]{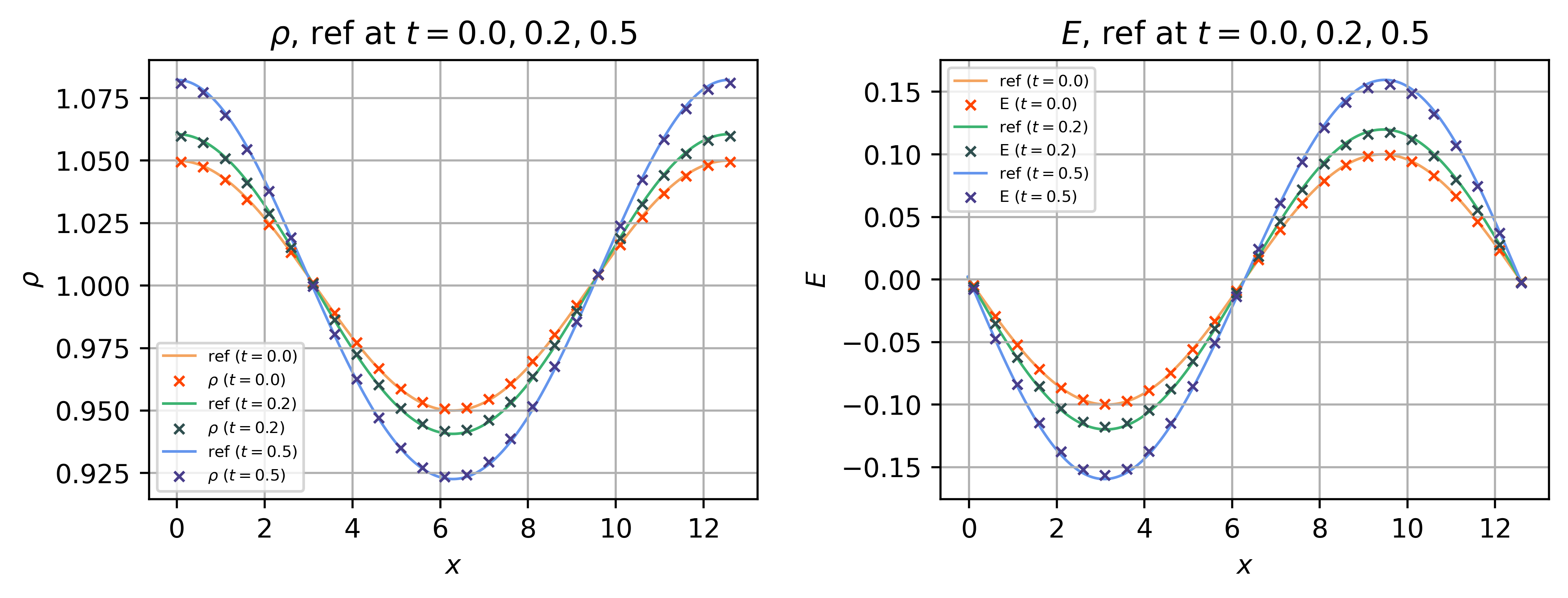}}
	\caption{The gravitational case solved by the micro-macro decomposition based and mass conservation based APNN methods with high-field regime ($\varepsilon=0.001$). Density $\rho$ (left column) and electric field $E$ (right column) as functions of space $x$ at $t = 0.0, 0.2, 0.5$. $(a)$ Neural networks are $[3, 128, 128, 128, 128, 128, 1]$ for $\rho, \phi$ and $[4, 256, 256, 256, 256, 256, 1]$ for $g$. Batch size is $512$ in domain and $256$ for initial condition.  Penalty $\lambda_1 = \lambda_3 = 1$. $(b)$ Neural networks are $[3, 128, 128, 128, 128, 128, 1]$ for $\rho, \phi$ and $[4, 256, 256, 256, 256, 256, 1]$ for $f$. Batch size is $512$ in domain, $256$ for initial condition and $256$ for conservation condition $\rho=\langle f\rangle$. Penalty $\kappa_1 = 50$; and $\kappa_3 = \kappa_4 = 1$. Errors: $(a)$ Relative $\ell^2$ error of $\rho$ is $2.72\times 10^{-4}$ at $t=0.2$ and $9.78\times 10^{-4}$ at $t=0.5$. Relative $\ell^2$ error of $E$ is $3.00\times 10^{-3}$ at $t=0.2$ and $1.77\times 10^{-2}$ at $t=0.5$. $(b)$ Relative $\ell^2$ error of $\rho$ is $5.82\times 10^{-4}$ at $t=0.2$ and $8.75\times 10^{-4}$ at $t=0.5$. Relative $\ell^2$ error of $E$ is $1.52\times 10^{-2}$ at $t=0.2$ and $1.94\times 10^{-2}$ at $t=0.5$.}
    \label{gravity}
\end{figure}

The initial and boundary conditions applied are the same as in the Landau damping case. However, $\phi_0$ is now the analytical solution to the equation $\partial_x\phi_0(x)= \rho_0(x)-h(x)$ due to the gravitational force. As shown in Figure \ref{gravity}, the solutions obtained from both APNN methods exhibit good agreement with the exact solution. We conclude that both APNN methods correctly capture the dynamics of the VPFP system with gravitational force.

%%%%%%%%%%%%%%%%%%%% UQ %%%%%%%%%%%%%%%%%%%%%
\subsection{Problem VI: Uncertainty quantification (UQ) problems}
This section highlights the potential of the proposed APNN methods for high-dimensional problems. Such problems present major challenges for classical mesh-based numerical schemes due to substantial computational costs. Here, we consider the uncertainty quantification (UQ) problems with high-dimensional random inputs. A typical setup of UQ is represented by the high-dimensional random input vectors. The initial conditions with random inputs are set as:
\begin{equation}
    \rho_0(x,\mathbf{z})=\frac{\sqrt{2 \pi}}{2}(2+\cos (2 \pi x))+0.1\prod_{i=1}^{10} \sin \left(\pi z_i\right), \,\,\, \boldsymbol{z}=\left(z_1, z_2, \cdots, z_{10}\right) \sim \mathscr{U}\left([-1,1]^{10}\right),
\end{equation}
where $z_1, z_2, \cdots, z_{10}$ are independent random variables following a uniform distribution over $[-1,1]$.
The background charge $h$ is set as:
\begin{equation}
    h(x,\mathbf{z})=\sqrt{2 \pi}+0.1\prod_{i=1}^{10} \sin \left(\pi z_i\right), \,\,\, \boldsymbol{z}=\left(z_1, z_2, \cdots, z_{10}\right) \sim \mathscr{U}\left([-1,1]^{10}\right).
\end{equation}
Then one analytically solves initial $\phi_0$ by the Poisson equation and defines:
\begin{equation}
    f_0(x, v, \mathbf{z})=\frac{\rho_0}{\sqrt{2 \pi}} \exp(-\frac{\left(v+\partial_x \phi_0\right)^2}{2}).
\end{equation}
The computational domain is $(x, v)\in[0,1]\times[-6, 6]$. Periodic boundary condition in $x$ direction is applied. The high-field limit characterized by $\varepsilon=0.001$ is considered.

\begin{figure}[htbp]
	\centering  %图片全局居中
	\subfigbottomskip=6pt %两行子图之间的行间距
	\subfigcapskip=-8pt %设置子图与子标题之间的距离
	\subfigure[Kinetic regime $(\varepsilon=1)$.]{
		\includegraphics[width=0.6\linewidth]{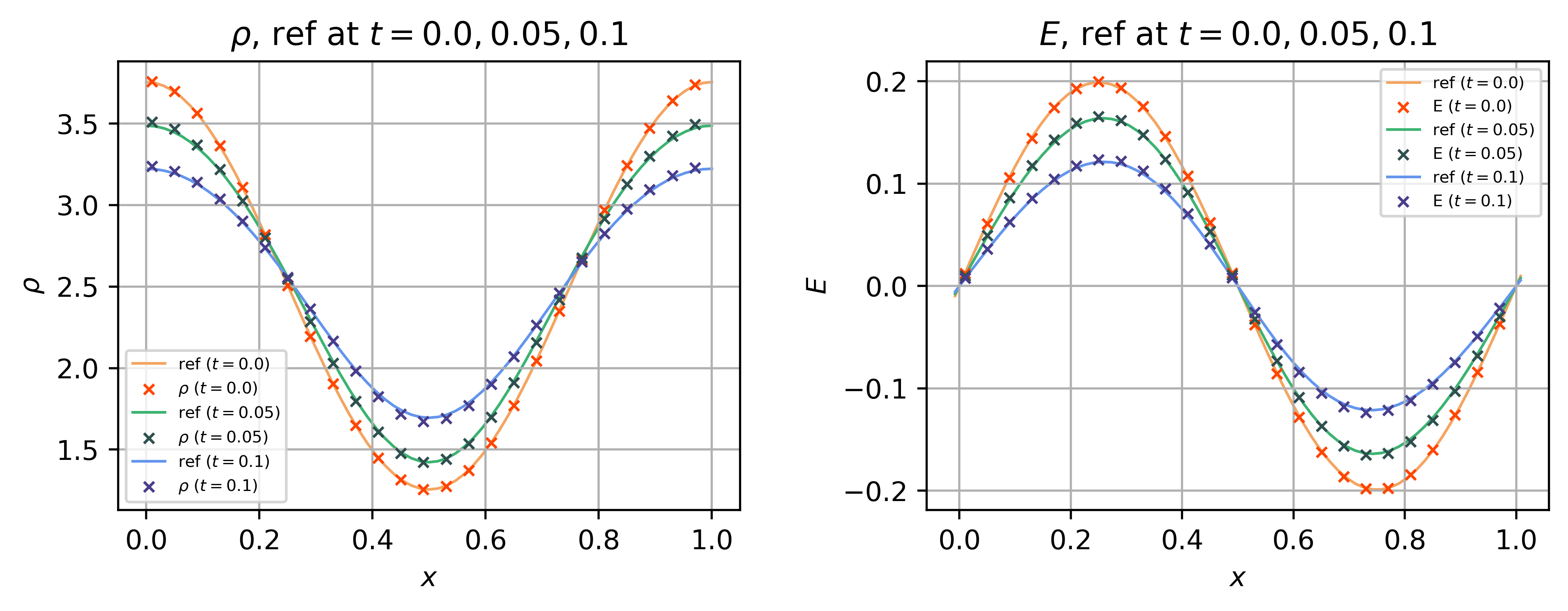}}
	\subfigure[High-field regime $(\varepsilon=0.001)$.]{
		\includegraphics[width=0.6\linewidth]{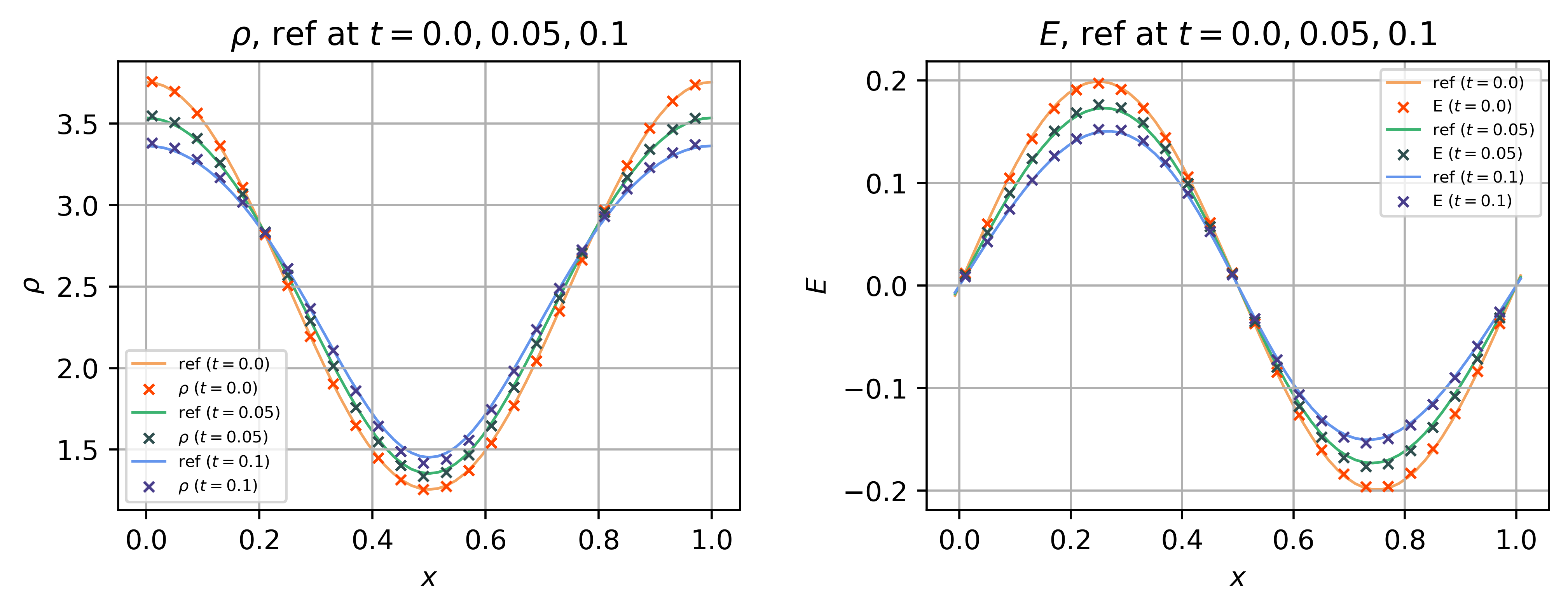}}
	\caption{UQ problem solved by the micro-macro decomposition based APNN method. Density $\rho$  (left column) and electric field $E$ (right column) as functions of space $x$ at $t = 0.0, 0.05, 0.1$. Neural networks are $[13, 128, 128, 128, 128, 128, 1]$ for $\rho, \phi$ and $[14, 256, 256, 256, 256, 256, 1]$ for $g$. Batch size is $512$ in domain and $256$ for initial condition. Penalty $\lambda_1 = 50$, $\lambda_3 = 1$ for kinetic regime; and $\lambda_1 = 1$, $\lambda_3 = 1$ for high-field regime. Errors: $(a)$ Relative $\ell^2$ error of $\rho$ is $4.71\times 10^{-3}$ at $t=0.05$ and $5.00\times 10^{-3}$ at $t=0.1$. Relative $\ell^2$ error of $E$ is $1.30\times 10^{-2}$ at $t=0.05$ and $2.28\times 10^{-2}$ at $t=0.1$. $(b)$ Relative $\ell^2$ error of $\rho$ is $4.67\times 10^{-3}$ at $t=0.05$ and $7.86\times 10^{-3}$ at $t=0.1$. Relative $\ell^2$ error of $E$ is $2.11\times 10^{-2}$ at $t=0.05$ and $1.77\times 10^{-2}$ at $t=0.1$.}
    \label{uq mm}
\end{figure}

\begin{figure}[htbp]
	\centering  %图片全局居中
	\subfigbottomskip=6pt %两行子图之间的行间距
	\subfigcapskip=-8pt %设置子图与子标题之间的距离
	\subfigure[Kinetic regime $(\varepsilon=1)$.]{
		\includegraphics[width=0.6\linewidth]{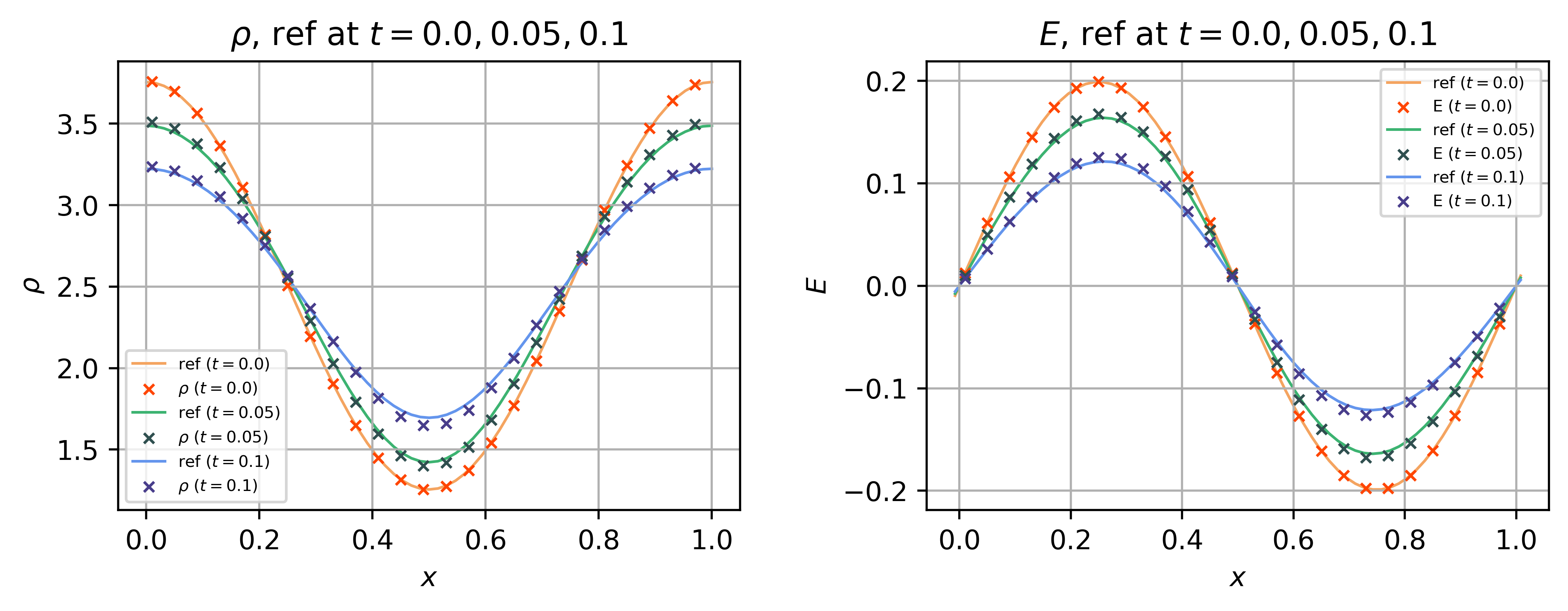}}
	\subfigure[High-field regime $(\varepsilon=0.001)$.]{
		\includegraphics[width=0.6\linewidth]{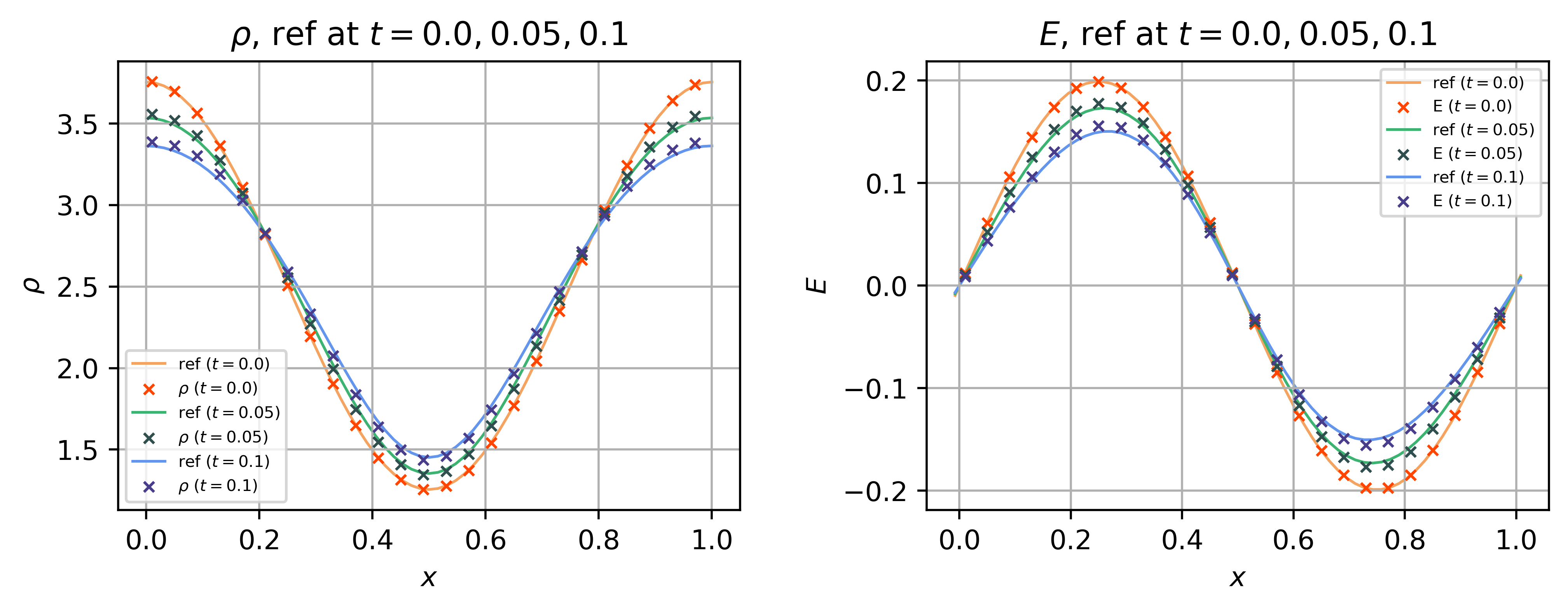}}
	\caption{UQ problem solved by the mass conservation based APNN method. Density $\rho$ (left column) and electric field $E$ (right column) as functions of space $x$ at $t = 0.0, 0.05, 0.1$. Neural networks are $[13, 128, 128, 128, 128, 128, 1]$ for $\rho, \phi$ and $[14, 256, 256, 256, 256, 256, 1]$ for $f$. Batch size is $512$ in domain, $256$ for initial condition and $256$ for conservation condition $\rho=\langle f\rangle$. Penalty $\kappa_1 = 100$, $\kappa_3 = \kappa_4 = 1$ for kinetic regime and $\kappa_1 = 100$, $\kappa_3 = \kappa_4 = 1$ for high-field regime. Errors: $(a)$ Relative $\ell^2$ error of $\rho$ is $7.83\times 10^{-3}$ at $t=0.05$ and $1.01\times 10^{-2}$ at $t=0.1$. Relative $\ell^2$ error of $E$ is $2.58\times 10^{-2}$ at $t=0.05$ and $4.02\times 10^{-2}$ at $t=0.1$. $(b)$ Relative $\ell^2$ error of $\rho$ is $7.57\times 10^{-3}$ at $t=0.05$ and $1.11\times 10^{-2}$ at $t=0.1$. Relative $\ell^2$ error of $E$ is $2.57\times 10^{-2}$ at $t=0.05$ and $3.34\times 10^{-2}$ at $t=0.1$.}
    \label{uq mc}
\end{figure}

We present the performance of the two proposed APNN methods in the high-demensional case. Density $\rho$ and electric field $E$ are evaluated at $t=0.0, 0.05, 0.1$ by taking expectation on $10^4$ times simulations for $(z_1, \cdots ,z_{10})$. Figure \ref{uq mm} and Figure \ref{uq mc} exhibit good approximations between the APNN approximations and the exact solutions. These agreements imply the potential for APNN methods to resolve high-dimensional and multiscale problems.

%%%%%%%%%%%%%%%%%%%%%%%%%%%%%%%%%%%%%%%%%%%%%%%%
%%%%%%%%%%%%%%%%%%%%%%%%%%%%%%%%%%%%%%%%%%%%%%%%
\section{Conclusion}

In this work, we developed two numerical methods based on PINN--with tailored loss function suitable for solving the multiscale uncertain Vlasov-Poisson-Fokker-Planck system with possible the high-field regimes. Our work incorporates the Asymptotic-Preserving (AP) mechanism into the loss function. Besides the micro-macro decomposition based Asymptotic-Preserving neural network (APNN) method, we formulate another APNN method that enforces the mass conservation law. This  method exhibits high adaptability for long time duration and non-equilibrium initial data. Diverse numerical examples are conducted across varied regimes to validate and compare the performance of both APNN methods. The numerical examples demonstrate the effectiveness of two proposed APNN  methods for approximating multiscale VPFP systems. While the micro-macro decomposition based APNN method shows its accuracy in discountinuity cases, the mass conservation based method shows its broader applicability for the cases with long-duration or non-equilibrium initial data.

\section*{Acknowledgement}
Shi Jin is partially supported by the Strategic Priority Research Program of Chinese Academy of Sciences XDA25010401, the NSFC grant No. 12031013, the Shanghai Municipal Science and Technology Major Project (2021SHZDZX0102), and the Fundamental Research Funds for the Central Universities. Zheng Ma is supported by NSFC Grant No. 12031013, No. 92270120 and Foundation of LCP.

\bibliographystyle{siam}
%\bibliography{vpfp}
\end{document}